\documentclass[a4paper,12pt]{amsart}

\usepackage{amsmath}
\usepackage{amssymb}
\usepackage{mathrsfs}
\usepackage{ifthen}
\usepackage{graphicx}
\usepackage[T1]{fontenc} 

\def\switchlinenumbers{\@ifstar
    {\let\makeLineNumberOdd\makeLineNumberRight
     \let\makeLineNumberEven\makeLineNumberLeft}%
    {\let\makeLineNumberOdd\makeLineNumberLeft
     \let\makeLineNumberEven\makeLineNumberRight}%
    }

\def\setmakelinenumbers#1{\@ifstar
  {\let\makeLineNumberRunning#1%
   \let\makeLineNumberOdd#1%
   \let\makeLineNumberEven#1}%
  {\ifx\c@linenumber\c@runninglinenumber
      \let\makeLineNumberRunning#1%
   \else
      \let\makeLineNumberOdd#1%
      \let\makeLineNumberEven#1%
   \fi}%
  }

\nonstopmode \numberwithin{equation}{section}
\newtheorem*{theorem*}{Theorem}

\newtheorem{thm}{Theorem}[section]
\newtheorem{cor}[equation]{Corollary}
\newtheorem{lem}[equation]{Lemma}
\newtheorem{prop}[equation]{Proposition}

\theoremstyle{definition}
\newtheorem{defn}{Definition}[section]
\newtheorem{example}{Example}[section]
\newtheorem{qsn}[equation]{Question}
\newtheorem{prob}[equation]{Problem}
\newtheorem{rem}{Remark}[section]

\newcounter{minutes}
\divide\time by 60
\newcounter{hours}
\multiply\time by 60
\addtocounter{minutes}{-\time}

\newcounter {own}
\def\theown {\thesection       .\arabic{own}}

\newenvironment{pf}[1][]{%
 \vskip 3mm
 \noindent
 \ifthenelse{\equal{#1}{}}%
  {{\slshape Proof. }}%
  {{\slshape #1.} }%
 }%
{\qed\bigskip}

\newcounter{alphabet}

\def\be{\begin{equation}}
\def\ee{\end{equation}}

\newcommand{\bee}{\begin{enumerate}}
\newcommand{\eee}{\end{enumerate}}

\newcommand{\blem}{\begin{lem}}
\newcommand{\elem}{\end{lem}}
\newcommand{\bthm}{\begin{thm}}
\newcommand{\ethm}{\end{thm}}
\newcommand{\bcor}{\begin{cor}}
\newcommand{\ecor}{\end{cor}}
\newcommand{\beg}{\begin{examp}}
\newcommand{\eeg}{\end{examp}}
\newcommand{\begs}{\begin{examples}}
\newcommand{\eegs}{\end{examples}}
\newcommand{\bdefe}{\begin{defin}}
\newcommand{\edefe}{\end{defin}}
\newcommand{\bprob}{\begin{prob}}
\newcommand{\eprob}{\end{prob}}
\newcommand{\bei}{\begin{itemize}}
\newcommand{\eei}{\end{itemize}}
\newcommand{\real}{{\operatorname{Re}\,}}

\begin{document}

\title{On Bloch norm and Bohr phenomenon for harmonic Bloch functions on simply connected domains}

\author{Vasudevarao Allu}
\address{Vasudevarao Allu,
School of Basic Sciences,
Indian Institute of Technology Bhubaneswar,
Bhubaneswar-752050, Odisha, India.}
\email{avrao@iitbbs.ac.in}

\author{Himadri Halder}
\address{Himadri Halder,
School of Basic Sciences,
Indian Institute of Technology Bhubaneswar,
Bhubaneswar-752050, Odisha, India.}
\email{himadrihalder119@gmail.com}

\subjclass[{AMS} Subject Classification:]{Primary 30A10, 30B10, 30C50, 30C55, 30C62, 30F45, 30H30, 31A05; Secondary 30C20, 30C65}

\def\thefootnote{}
\footnotetext{ {\tiny File:~\jobname.tex,
printed: \number\year-\number\month-\number\day,
          \thehours.\ifnum\theminutes<10{0}\fi\theminutes }
} \makeatletter\def\thefootnote{\@arabic\c@footnote}\makeatother

\begin{abstract}
For $\alpha \in (0,\infty)$, let $\mathcal{B}_{\mathcal{H},\Omega}(\alpha)$ denote the class of $\alpha$-Bloch mappings on a proper simply connected domain $\Omega \subseteq \mathbb{C}$. 
In this article, we introduce the class $\mathcal{B}^{*}_{\mathcal{H},\Omega}(\alpha)$ of harmonic $\alpha$-Bloch-type mappings on a proper simply connected domain $\Omega \subseteq \mathbb{C}$ and study several interesting properties of the classes $\mathcal{B}_{\mathcal{H},\Omega}(\alpha)$ and $\mathcal{B}^{*}_{\mathcal{H},\Omega}(\alpha)$ when $\Omega$ is proper simply connected domain and the shifted disk $\Omega_{\gamma}$ containing $\mathbb{D}$, where 
	$$
	\Omega_{\gamma}:=\bigg\{z\in\mathbb{C} : \bigg|z+\frac{\gamma}{1-\gamma}\bigg|<\frac{1}{1-\gamma}\bigg\}
	$$
	and $0 \leq \gamma <1$. We establish the Landau's theorem for the harmonic Bloch space $\mathcal{B}_{\mathcal{H},\Omega _{\gamma}}(\alpha)$ on the shifted disk $\Omega_{\gamma}$. 
	 For $f \in \mathcal{B}_{\mathcal{H},\Omega}(\alpha)$ (respectively $\mathcal{B}^{*}_{\mathcal{H},\Omega}(\alpha)$) of the form $f(z)=h(z) + \overline{g(z)}=\sum_{n=0}^{\infty}a_nz^n +  \overline{\sum_{n=1}^{\infty}b_nz^n}$ in $\mathbb{D}$ with Bloch norm $||f||_{\mathcal{H},\Omega, \alpha} \leq 1$ (respectively $||f||^{*}_{\mathcal{H},\Omega, \alpha} \leq 1$), we define the Bloch-Bohr radius for the class $\mathcal{B}_{\mathcal{H},\Omega}(\alpha)$ (respectively $\mathcal{B}^{*}_{\mathcal{H},\Omega}(\alpha)$) to be the largest radius $r_{\Omega,f} \in (0,1)$ such that $\sum_{n=0}^{\infty}(|a_n|+|b_{n}|) r^n\leq 1$ for $r \leq r_{\Omega, \alpha}$ and for all $f \in \mathcal{B}_{\mathcal{H},\Omega}(\alpha)$ (respectively $\mathcal{B}^{*}_{\mathcal{H},\Omega}(\alpha)$). We also investigate Bloch-Bohr radius for the classes $\mathcal{B}_{\mathcal{H},\Omega}(\alpha)$ and $\mathcal{B}^{*}_{\mathcal{H},\Omega}(\alpha)$ on simply connected domain $\Omega$ containing $\mathbb{D}$. 
\end{abstract}

\maketitle
\pagestyle{myheadings}
\markboth{Vasudevarao Allu and  Himadri Halder}{On Bloch norm and Bohr phenomenon for harmonic Bloch functions on simply connected domains}

\section{Introduction}
 The Bloch spaces, pre-schwarzian and schwarzian norm, Bohr phenomenon and their various generalizations have become a central object of study and several outstanding problems remain unsolved in classical geometric function theory on $\mathbb{D}:=\{z \in \mathbb{C}: |z|<1\}$. For $\alpha \in (0, \infty)$, an analytic function $f$ in the unit disk $\mathbb{D}$ is called a $\alpha$-Bloch function if 
$$
\beta (\alpha):=\sup _{z \in \mathbb{D}} (1-|z|^2)^\alpha |f'(z)|< \infty.
$$
The class of all $\alpha$-Bloch functions is denoted by $\mathcal{B}(\alpha)$. The class $\mathcal{B}(\alpha)$ is a Banach space with respect to the Bloch norm $|| .\,\, ||_{ \alpha}$, which is defined by $|| f\, ||_{\alpha}= |f(0)|+\beta (\alpha)$. The $\alpha$-Bloch space is a generalization of classical Bloch space $\mathcal{B}(1)$.
\vspace{3mm}

The definition of $\alpha$-Bloch space can be generalized to an arbitrary proper simply connected domain $\varOmega$ in $\mathbb{C}$ by means of hyperbolic metric. The hyperbolic metric $\lambda _{\mathbb{D}} |dz|$ of $\mathbb{D}$ is defined by $\lambda _{\mathbb{D}}(z)=1/(1-|z|^2)$, $z \in \mathbb{D}$. The equivalent form of the norm $\beta(\alpha)$ for $(0,\infty)$ in terms of hyperbolic density becomes 
$$
 \beta(\alpha)= \sup _{z \in \mathbb{D}} \frac{|f'(z)|}{\lambda ^{\alpha}_{\mathbb{D}}(z)}.
$$
Let $\Omega \subseteq \mathbb{C}$ be a proper simply connected domain and $f:\Omega \rightarrow \mathbb{D}$ be a conformal map. Then the hyperbolic metric $\lambda _{\Omega}(z) |dz|$ of $\Omega$ is defined by (see \cite{beardon-minda-hyperbolic-density})
\begin{equation} \label{h-v-p5-e-1.1}
\lambda _{\Omega}(z):=\lambda _{\mathbb{D}}(f(z)) |f'(z)|, \,\, z \in \Omega.
\end{equation}
We note that the definition \eqref{h-v-p5-e-1.1} is independent of the choice of conformal mapping of $\Omega$ onto $\mathbb{D}$. 
Let $f$ be an analytic function in $\varOmega \subseteq \mathbb{C}$. Then $f$ is said to be an $\alpha$-Bloch function in $\Omega$ if 
$$
\beta _{\Omega}(\alpha):=\sup_{z \in \Omega} \frac{|f'(z)|}{\lambda ^{\alpha}_{\Omega}(z)}<\infty
$$
and the space of all $\alpha$-Bloch functions is denoted by $\mathcal{B}_{\Omega}(\alpha)$. We define the norm on $\mathcal{B}_{\Omega}(\alpha)$ by $|| f\, ||_{\Omega, \alpha}= |f(0)|+\beta _{\Omega}(\alpha)$. In $1974$, Anderson {\it et al.} \cite{Anderson-1974} established several results on the coefficients and zeros of Bloch functions and the boundary behaviour of normal functions. In $1993$, Rohde \cite{Rohde-1993} studied the boundary behaviour of Bloch functions. In $1994$, Bonk {\it et al.} \cite{Bonk-CMFT-1994} studied extensively the hyperbolic metric on Bloch regions. There has been a significant work on the sharp distortion estimates for locally univalent Bloch functions in \cite{Bonk-J-Analyze-Math-1996,Yanigaha-BLMS-1994} and references therein. M\"{o}bius invariant space in the contest of Bloch spaces has been extensively studied by Arazy {\it et al.} in \cite{Arazy-1985}. Gnuschke-Hauschild and Pommerenke \cite{Gnuschke-Hauschild-1986} have established several interesting results on Bloch functions for gap series.
\vspace{3mm}

A complex-valued twice continuously differentiable function $f$ defined in a simply connected domain $\Omega \subseteq \mathbb{C}$ is said to be harmonic in $\Omega$ if $\triangle f=4f_{z\bar{z}}=0$ for all $z \in \Omega$, where $\triangle = \partial^{2}/\partial x^2 + \partial^{2}/\partial y^2 $. It is well-known that every harmonic mapping $f$ in $\Omega$ has the {\it canonical decomposition} $f=h+\overline{g}$, where $h$ and $g$ are analytic functions in $\Omega$ with $g(0)=0$. The Jacobian of $f$ is defined by $J_{f}=|h'|^2 - |g'|^2$. The function $f=h+\overline{g}$ is locally univalent and sense-preserving in $\Omega$ if, and only if, $J_{f}>0$ in $\Omega$ {\it i.e.,} $|h'|>|g'|$ or $|\omega _{f}|<1$ in $\Omega$, where $\omega_{f}=g'/h'$ is the dilation of $f$.
\vspace{3mm}

For a given $\alpha \in (0,\infty)$, a harmonic mapping $f=h+\overline{g}$ in the unit disk $\mathbb{D}$ is called an $\alpha$-Bloch mapping if 
$$
\beta _{\mathcal{H}}(\alpha):=\sup_{z \in \mathbb{D}} (1-|z|^2)^\alpha (|h'(z)|+|g'(z)|) < \infty.
$$
Then $\beta _{\mathcal{H}}(\alpha)$ defines a semi-norm, and the space equipped with the norm $ ||f\, ||_{\mathcal{H}, \alpha}= |f(0)|+\beta _{\mathcal{H}}(\alpha)$ is called the harmonic $\alpha$-Bloch space, denoted by $\mathcal{B}_{\mathcal{H}}(\alpha)$.  The space $\mathcal{B}_{\mathcal{H}}(\alpha)$ is a Banach space with respect to the norm $||.\, ||_{\mathcal{H},\alpha}$. In particular, when $g\equiv 0$, the space $\mathcal{B}_{\mathcal{H}}(\alpha)$ coincides with $\mathcal{B}(\alpha)$. Thus, $\mathcal{B}_{\mathcal{H}}(\alpha)$ is a generalization of $\mathcal{B}(\alpha)$. In $1989$, Colona \cite{colona-1989} studied extensively the space $\mathcal{B}_{\mathcal{H}}(1)$, which is a generalization of the classical Bloch space $\mathcal{B}(1)$.  The harmonic $\alpha$-Bloch space $\mathcal{B}_{\mathcal{H}}(\alpha)$ is a generalization of $\mathcal{B}_{\mathcal{H}}(1)$ (see \cite{chen-2011}). Motivated by the well-known results on analytic Bloch space, in $2016$, Efraimidis {\it et al.}\cite{efraimidis-2017} introduced harmonic $1$-Bloch-type mappings. For a given $\alpha \in (0,\infty)$, a harmonic mapping $f$ in $\mathbb{D}$ is called a harmonic $\alpha$-Bloch-type mapping if 
$$
\beta ^*_{\mathcal{H}}(\alpha):= \sup_{z \in \mathbb{D}} (1-|z|^2)^\alpha \sqrt{|J_{f}(z)|}<\infty,
$$
where $J_{f}$ is the Jacobian of $f$ defined by $J_{f}(z)=|h'(z)|^2 - |g'(z)|^2$. 
Let $\mathcal{B}^*_{\mathcal{H}}(\alpha)$ be the space of all $\alpha$-Bloch-type mappings and $||f\, ||^*_{\mathcal{H}, \alpha}:=|f(0)|+ \beta ^*_{\mathcal{H}}(\alpha)$ be the pseudo-norm of $f$.
\vspace{3mm}

Motivated by the $\alpha$-Bloch mappings and $\alpha$-Bloch-type mappings in $\mathbb{D}$, in this paper, we consider  $\alpha$-Bloch mappings and $\alpha$-Bloch-type mappings in an arbitrary simply connected domain $\Omega \subseteq \mathbb{C}$. 
\begin{defn}
For a given $\alpha \in (0,\infty)$, a harmonic mapping $f$ in a simply connected domain $\Omega$ is called a harmonic $\alpha$-Bloch mapping if 
$$
\beta_{\mathcal{H}, \Omega} (\alpha):=\sup_{z \in \Omega} \frac{|h'(z)|+|g'(z)|}{\lambda^{\alpha}_{\Omega}(z)}< \infty.
$$
\end{defn}
We define the class of all $\alpha$-Bloch mappings by $\mathcal{B}_{\mathcal{H},\Omega}(\alpha)$ and the Bloch norm is defined by $||f\, ||_{\mathcal{H}, \Omega, \alpha}:=|f(0)|+ \beta _{\mathcal{H}, \Omega}(\alpha)$.
\begin{defn}
	 For a given $\alpha \in (0,\infty)$, a harmonic mapping $f$ in a simply connected domain $\Omega$ is called a harmonic $\alpha$-Bloch-type mapping if 
	$$
	\beta^{*}_{\mathcal{H}, \Omega} (\alpha):=\sup_{z \in \varOmega} \frac{\sqrt{|J_{f}(z)|}}{\lambda^{\alpha}_{\Omega}(z)}< \infty.
	$$
\end{defn}
We define the class of all $\alpha$-Bloch-type mappings by $\mathcal{B}^*_{\mathcal{H},\Omega}(\alpha)$ and we denote the pseudo-norm by $||f\, ||^*_{\mathcal{H}, \Omega, \alpha}:=|f(0)|+ \beta^* _{\mathcal{H}, \Omega}(\alpha)$. Observe that for $\Omega=\mathbb{D}$,  $\mathcal{B}_{\mathcal{H},\Omega}(\alpha)$ and $\mathcal{B}^*_{\mathcal{H},\Omega}(\alpha)$ coincide with the spaces $\mathcal{B}_{\mathcal{H}}(\alpha)$ and $\mathcal{B}^*_{\mathcal{H}}(\alpha)$ respectively. It is important to note that 
$$
\frac{\sqrt{|J_{f}(z)|}}{\lambda^{\alpha}_{\Omega} (z)} \leq \frac{|h'(z)|+|g'(z)|}{\lambda^{\alpha}_{\Omega}(z)} \,\,\,\,\, \mbox{for} \,\, z\in \mathbb{D},
$$
which clearly shows that $\mathcal{B}_{\mathcal{H},\Omega}(\alpha) \subseteq \mathcal{B}^*_{\mathcal{H},\Omega}(\alpha)$ and hence, $\mathcal{B}^*_{\mathcal{H},\Omega}(\alpha)$ is a generalization of $\mathcal{B}_{\mathcal{H},\Omega}(\alpha)$. In particular, when $g \equiv 0$ {\it i.e.}, $f$ is analytic function in $\Omega$, the spaces $\mathcal{B}_{\mathcal{H},\Omega}(\alpha)$ and $\mathcal{B}^*_{\mathcal{H},\Omega}(\alpha)$ coincide with analytic Bloch space $\mathcal{B}_{\Omega}(\alpha)$. Therefore, we have 
$$
||f\, ||_{\mathcal{H}, \Omega, \alpha}=||f\, ||^*_{\mathcal{H}, \varOmega, \alpha}=||f\, ||_{\Omega, \alpha}.
$$

Bohr's famous power series theorem (see \cite{Bohr-1914}) asserts that if $f$ is an analytic function in $\mathbb{D}$ with the power series expansion $f(z)=\sum_{n=0}^{\infty} a_{n}z_{n}$ in $\mathbb{D}$ such that $|f(z)|<1$, then 
\begin{equation} \label{h-v-p5-e-1.2}
\sum \limits_{n=0}^{\infty} |a_{n}| |z|^n \leq 1 \,\, \mbox{for} \,\, |z|=r \leq 1/3.
\end{equation}
M. Riesz, Schur, and Weiner (see \cite{paulsen-2002}) have independently shown that $1/3$ is the best possible constant. The inequality \eqref{h-v-p5-e-1.2} is usually known as Bohr inequality and the constant $1/3$ is the Bohr radius for the class of analytic functions whose modulus is less than $1$. 
\vspace{3mm}

In $1995$, Dixon \cite{Dixon & BLMS & 1995} used Bohr inequality in connection with the long-standing open problem of characterizing Banach algebras satisfying the von Neumann inequality. In $1997$, Boas and Khavinson \cite{boas-1997} extended the Bohr inequality \eqref{h-v-p5-e-1.2} to several complex variables. Indeed, Boas and Khavinson \cite{boas-1997} have introduced the Bohr radius $K_{n}$ for the Hardy space $H^{\infty}(\mathbb{D}^n)$ of bounded holomorphic functions on the $n$-dimensional polydisc and proved that, if $n>1$, then 

\begin{equation} \label{h-v-p5-e-1.3}
\frac{1}{3\sqrt{n}}<K_{n}< \frac{2\sqrt{log\,n}}{\sqrt{n}}.
\end{equation}
 
The article \cite{boas-1997} by Boas and Khavinson, is a source of inspiration for many subsequent papers, connecting the asymptotic behaviour of $K_{n}$ to various problems in functional analysis, for example, geometry of Banach spaces, unconditional basis constant of spaces of polynomials. Thus, in the recent years, there has been a big interest in determining the behaviour of $K_{n}$ for the large values of $n$. In $2006$, Defant and Frerick \cite{defant-2006} improved the left inequality of \eqref{h-v-p5-e-1.3} to 
$K_{n} \geq c \sqrt{log \, n/(nlog\,log\,n)}$. Using the hypercontractivity of the polynomial Bohnenblust-Hille inequality, Defant {\it et al.}\cite{defant-2011} have proved that 
$$
K_{n}=b_{n} \sqrt{\frac{log\,n}{n}} \,\,\,\,\, \mbox{with}\, \, \frac{1}{\sqrt{2}}+ o(1) \leq b_{n} \leq 2.
$$
In $2014$, Bayart {\it et al.} \cite{bayart-advance-2014} obtained the exact asymptotic behaviour of $K_{n}$ and showed that 
$$
K_{n} \sim + \sqrt{\frac{log\,n}{n}} \, \, \mbox{or} \,\, \lim\limits_{n\rightarrow \infty} \frac{K_{n}}{\sqrt{\frac{log\,n}{n}}}=1.
$$
In $2019$, Popescu \cite{popescu-2019} extended the inequality \eqref{h-v-p5-e-1.3} for free holomorphic functions to polyballs $B_{n}$, $n=(n_{1}, n_{2}, \ldots, n_{k}) \in \mathbb{N}^k$, which is a non-commutative anlogue of the scalar polyball $(\mathbb{C}^{n_{1}} \times \cdots \times \mathbb{C}^{n_{k}}) \in \mathbb{N}^k $ and showed that 
$$
\frac{1}{3\sqrt{k}}<K_{m,h}(B_{n})< \frac{2\sqrt{log\,k}}{\sqrt{k}} \,\, \mbox{for} \,\, k>1,
$$
where $K_{m,h}(B_{n})$ is the Bohr radius associated with multi-homogeneous power series expansion of the free holomorphic functions. Aizenberg {\it et al.}\cite{aizenberg-2001} have generalized the Bohr theorem for bases in spaces of holomorphic functions of several complex variables. Bohr power series theorem connection with local Banach space theory has been extensively studied  by Defant {\it et al.} in \cite{defant-2003}. In $2011$, Defant {\it et al.}\cite{defant-JRAM-2011} have estimated Bohr radius in the unit ball of $l^{n}_{p}$. Recently, Liu and Ponnusamy \cite{Liu-Pon-PAMS-2020} have obtained several multidimensional analogues of refined Bohr inequality. For further work on multidimensional Bohr work, we refer \cite{aizn-2000,aizn-2007}. In $2021$, Bhowmik and Das \cite{bhowmik-2021} extensively studied Bohr inequalities for operator valued functions, which can be viewed as the analogues of a couple of interesting results from scalar valued settings. For more intriguing aspects of Bohr inequality, we refer \cite{abu-2011,Ali-2017,alkhaleefah-2019,Himadri-Vasu-P2,Himadri-Vasu-P3,Ayt & Dja & BLMS & 2013,bene-2004,Ismagilov-2020,Kayumov-Ponnusamy-2018-b} and references therein.
\vspace{3mm}

The Bohr phenomenon \cite{Abu-2010} for harmonic functions  $f$ of the form  $f(z)=h(z)+\overline {g(z)}$ in $\mathbb{D}$, where $h(z)=\sum_{n=0}^{\infty} a_{n}z^{n}$ 
and $g(z)=\sum_{n=1}^{\infty} b_{n}z^{n}$ is to find the largest radius $r_{f}$, $0<r_{f}<1$ such that 
\begin{equation} \label{h-v-p5-e-1.4}
\sum\limits_{n=1}^{\infty} (|a_{n}|+|b_{n}|) |z|^{n}\leq 1 
\end{equation}
holds for $|z|\leq r_{f}$. Kayumov {\it et al.} \cite{Kayumov-Ponnuswamy-MN-2018} have generalized the Bohr inequality for locally univalent harmonic mappings in $\mathbb{D}$. In $2018$, Liu and Ponnusamy \cite{Liu-Results-Math-2018} studied the Bohr inequality for harmonic $\nu$-Bloch and $\nu$-Bloch-type mappings in the unit disk $\mathbb{D}$.
Bohr phenomenon for certain subclasses of harmonic mappings has been studied in \cite{Ahamed-Allu-Halder-AAMP-2020,Himadri-Vasu-P1}. For more work on Bohr inequality for harmonic mappings, we refer \cite{abu-2014,Kayumov-Ponnusamy-2018-b,Liu-Ponnusamy-BMMS-2019}
\vspace{3mm}

In $2010$, Fournier and Ruscheweyh \cite{Four-Rusc-2010} estimated the Bohr radius in an arbitrary simply connected domain containing $\mathbb{D}$. Let $\mathcal{H}(\Omega)$ be the class of analytic functions in $\Omega$ and let $\mathcal{B}(\Omega)$ denote the class of functions $f \in \mathcal{H}(\Omega) $ such that $f(\Omega) \subseteq \overline{\mathbb{D}}$. The Bohr radius $B_{\Omega}$ for the class $\mathcal{B}(\Omega)$ is defined by 
$$
B_{\Omega}:=\sup\bigg\{r\in (0,1) : M_{f}(r)\leq 1\; \text{for all}\; f(z)=\sum_{n=0}^{\infty}a_nz^n\in\mathcal{B}(\Omega),\; z\in\mathbb{D}\bigg\},
$$
where $M_{f}(r):=\sum_{n=0}^{\infty}|a_n|r^n$ is the associated majorant series of $f \in \mathcal{B}(\Omega)$ in $\mathbb{D}$. Clearly, when $\Omega=\mathbb{D}$, $B_{\mathbb{D}}=1/3$, which is the classical Bohr radius for the class $\mathcal{B}(\mathbb{D})$.
\vspace{2mm}

For analytic functions $f$ in a simply connected domain $\varOmega$ containing $\mathbb{D}$, Fournier and  Ruscheweyh  \cite{Four-Rusc-2010} have estimated the Bohr radius $B_{\varOmega}$. For brevity, we recall the following result.
\begin{thm}\cite{Four-Rusc-2010}
	Let $ \Omega $ be a simply connected domain which contains the unit disk $ \mathbb{D} $ and let 
	\begin{equation}\label{e-2.3}
	\lambda:=\lambda(\Omega)=\sup_{f\in\mathcal{B}(\Omega),\;n\geq 1}\bigg\{\frac{|a_n|}{1-|a_0|^2} :  a_0 \not \equiv f(z)=\sum_{n=0}^{\infty}a_nz^n\;\; \mbox{for}\;\; z\in\mathbb{D}\bigg\}.
	\end{equation}
	Then $ 1/(1+2\lambda)\leq\mathcal{B}_{\Omega} $ and the equality $ \sum_{n=0}^{\infty}|a_n|(1/(1+2\lambda))^n=1 $ holds for a function $ f(z)=\sum_{n=0}^{\infty}a_nz^n $ in $ \mathcal{B}(\Omega) $ if, and only if, $ f\equiv c $ with $ |c|=1 $.
\end{thm}
In particular, when $\Omega=\Omega _{\gamma}$ is a disk containing $\mathbb{D}$ defined by 
$$
\Omega_{\gamma}:=\bigg\{z\in\mathbb{C} : \bigg|z+\frac{\gamma}{1-\gamma}\bigg|<\frac{1}{1-\gamma}\bigg\},
$$
Fournier and  Ruscheweyh  \cite{Four-Rusc-2010} have obtained the exact Bohr radius for the class $\mathcal{B}(\Omega_{\gamma})$.
\begin{thm}\cite{Four-Rusc-2010} \label{thm-1.2}
	For $ 0\leq \gamma<1 $, let $ f\in\mathcal{B}(\Omega_{\gamma}) $, with $ f(z)=\sum_{n=0}^{\infty}a_nz^n $ in $ \mathbb{D} $. Then,
	\begin{equation*}
	\sum_{n=0}^{\infty}|a_n|r^n\leq 1\;\; \text{for}\;\; r\leq\rho _{\gamma}:=\frac{1+\gamma}{3+\gamma}.
	\end{equation*}
	Moreover, $ \sum_{n=0}^{\infty}|a_n|\rho _{\gamma}^n=1 $ holds for a function $ f(z)=\sum_{n=0}^{\infty}a_nz^n $ in $ \mathcal{B}(\Omega_{\gamma}) $ if, and only if, $ f(z)=c $ with $ |c|=1 $.
\end{thm} 
Recently, Evdoridis {\it et al.} \cite{Evd-Ponn-Rasi-2020} have studied several improved version of Bohr inequality in the shifted disk $\Omega_{\gamma}$. Later, Ahamed {\it et al.} \cite{Ahamed-Allu-Halder-P3-2020} have established several improved version of Bohr radius, Bohr-Rogosinski radius and refined Bohr radius for the functions defined in $\Omega_{\gamma}$, and obtained several sharp results.  
\begin{rem}
\begin{enumerate}
\item
The sharp coefficients estimate play a vital role to obtain the sharp Bohr radius. The proofs of the inequality \eqref{h-v-p5-e-1.2} relied on the sharp coefficient inequalities which may be obtained as an application of Pick's invariant form of Schwarz's lemma for $ f\in\mathcal{B}(\mathbb{D}) $:
\begin{equation*}
|f^{\prime}(z)|\leq\frac{1-|f(z)|^2}{1-|z|^2}\;\;\mbox{for}\;\; z\in\mathbb{D}.
\end{equation*}  
In particular, $ |f^{\prime}(0)|=|a_1|\leq 1-|f(0)|^2=1-|a_0|^2 $ and hence from this, the sharp inequality $ |a_n|\leq 1-|a_0|^2 $ follows for $ n\geq 1 $. On the other hand, the proof of Theorem \ref{thm-1.2} follows from the sharp coefficients estimate $|a_{n}|\leq (1-|a_{0}|^2)/(1+\gamma)$ (see \cite{Four-Rusc-2010}). 

\item   It is pertinent to note that there is no assurance that the Bohr radius exists for every class of functions. For instance, Aizenberg \cite{aizenberg-2012} has shown that Bohr radius does not exist for the space of analytic functions defined in the annulus $\{z \in \mathbb{C}: t<|z|<1, \, 0<t<1\}$.
\end{enumerate}

\end{rem}
The sharp coefficient bounds for the Bloch spaces $\mathcal{B}_{\mathcal{H}}(\alpha), \mathcal{B}^*_{\mathcal{H}}(\alpha), \mathcal{B}_{\mathcal{H},\Omega}(\alpha) $, and $\mathcal{B}^*_{\mathcal{H},\Omega}(\alpha)$ are not yet known. Thus, the above two discussion led to the the following question.
\begin{qsn}
	Let the function $f=h+ \overline{g}$ belongs $\mathcal{B}_{\mathcal{H},\Omega}(\alpha)$ (respectively $\mathcal{B}^*_{\mathcal{H},\Omega}(\alpha)$)  with the restriction to $\mathbb{D}$ {\it i.e.,} $h(z)=\sum_{n=0}^{\infty} a_{n}z^{n}$ 
	and $g(z)=\sum_{n=1}^{\infty} b_{n}z^{n}$ in $\mathbb{D}$  such that $||f\, ||_{\mathcal{H}, \varOmega,\alpha} \leq 1$ ($||f\, ||^*_{\mathcal{H}, \varOmega,\alpha} \leq 1$). Can we establish the inequality \eqref{h-v-p5-e-1.4} for $|z|=r \le r_{f} \in (0,1)$, for all functions $f \in \mathcal{B}_{\mathcal{H},\Omega}(\alpha)$ (respectively $\mathcal{B}^*_{\mathcal{H},\Omega}(\alpha)$) without knowing the upper bounds for $|a_{n}|$ and $|b_{n}|$?
\end{qsn}
The answer to Question $1.6$ is affirmative, which we discuss in Section $4$. When such radius $r_{f}$ exists for $\mathcal{B}_{\mathcal{H},\Omega}(\alpha)$ (respectively $\mathcal{B}^*_{\mathcal{H},\Omega}(\alpha)$), we call it as Bloch-Bohr radius for class $\mathcal{B}_{\mathcal{H},\Omega}(\alpha)$ (respectively $\mathcal{B}^*_{\mathcal{H},\Omega}(\alpha)$).

\begin{figure}[!htb]
	\begin{center}
		\includegraphics[width=0.48\linewidth]{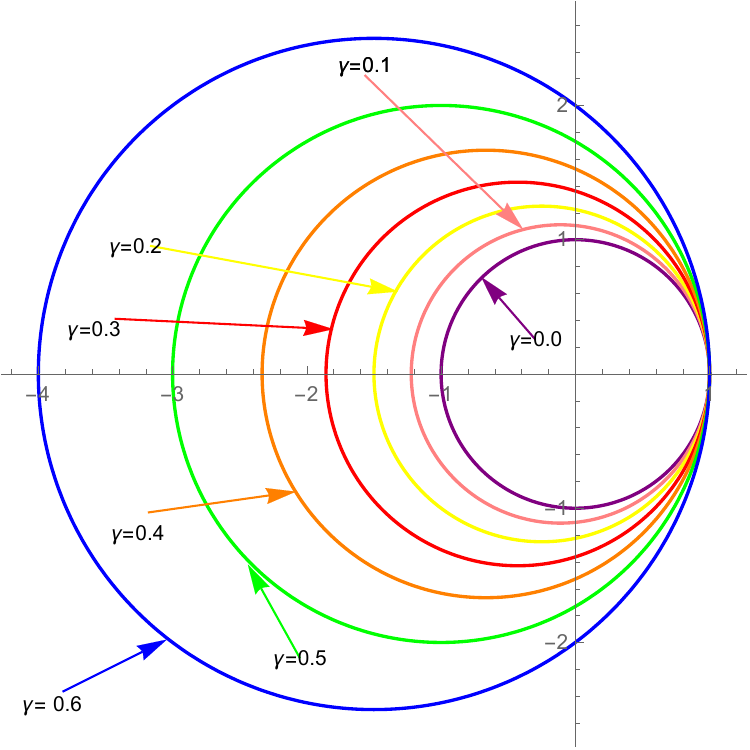}
	\end{center}
	\caption{The graph of $ C_{\gamma} $ when $\gamma=0,0.1, 0.2, 0.3,0.4, 0.5, 0.6$.}
	\label{figure-sec-1}
\end{figure}

It is important to note that $\Omega_{\gamma}$ contains $\mathbb{D}$ and is increasing in $\gamma \in [0,1)$. Here $\Omega_{\gamma}$ is increasing in the sense that if $\gamma_{1}, \gamma_{2} \in [0,1)$ such that $\gamma_{1} \leq \gamma_{2}$, then $\Omega _{\gamma_{1}} \subseteq \Omega _{\gamma_{2}}$. Figure \ref{figure-sec-1} shows that the pictures of the circles $C_{\gamma}: |z+\gamma/ (1-\gamma)|=1/(1-\gamma) $ for certain values of $\gamma \in [0,1)$.

From Theorem \ref{thm-1.2}, we have Bohr radius $\rho_{\gamma}=(1+\gamma)/(3+\gamma)$ for the restriction map $f$ in $\mathbb{D}$, when $f \in \mathcal{B}(\Omega_{\gamma})$. We calculate Bohr radius $\rho_{\gamma}$ for certain values of $\gamma \in [0,1)$ in Table \ref{tabel-sec-1}.

\begin{table}[ht]
	\centering
	\begin{tabular}{|l|l|l|l|l|l|l|l|l|l|l|}
		\hline
		$\gamma$& $ 0.0$& $0.1$ & $0.2$ & $0.3$ &$0.4$ &$0.5$ &$0.6$ &$0.7$ &$0.8$ &$0.9$ \\
		\hline
		$\rho_{\gamma}$& $ 0.3333$& $0.3548$ & $0.3750$ & $0.3939$ &$0.4117$ &$0.4285$ &$0.4444$ &$0.4594$ &$0.4736$ &$0.4871$\\
		\hline		
	\end{tabular}
	\vspace{3mm}
	\caption{Bohr radius $\rho_{\gamma}$ for the class $\mathcal{B}(\Omega_{\gamma})$ for different values of $\gamma \in [0,1)$.}
	\label{tabel-sec-1}
\end{table}

\begin{figure}[!htb]
	\begin{center}
		\includegraphics[width=0.40\linewidth]{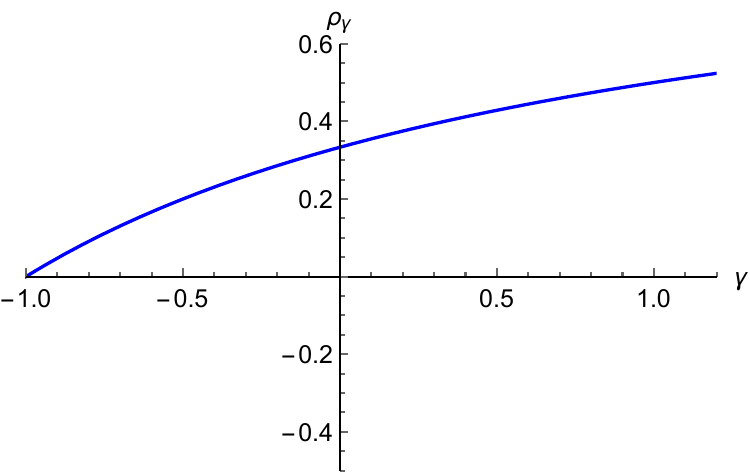}
	\end{center}
	\caption{The graph of $\rho_{\gamma} $ for $\gamma \in [0,1)$.}
	\label{figure-sec-1-a}
\end{figure}

From Table \ref{tabel-sec-1} and Figure \ref{figure-sec-1-a}, we observe that $\rho_{\gamma}$ is increasing in $\gamma \in [0,1)$. Indeed, $\rho'_{\gamma}= 2/(3+\gamma)^2 >0$ for $\gamma \in [0,1)$ and hence $\rho_{\gamma}$ is increasing in $\gamma \in [0,1)$. We note that when $\gamma_{1} \leq \gamma_{2}$, $\Omega_{\gamma_{1}} \subseteq \Omega_{\gamma_{2}}$, that is,  $\Omega_{\gamma}$ is increasing when $\gamma$ is increasing. Therefore, when $\Omega_{\gamma}$ is increasing then Bohr radius is also increasing accordingly. This fact leads us to consider the following questions:

\begin{qsn} \label{him-vasu-P5-qsn-1.7}
If the Bloch-Bohr radius exists for $\mathcal{B}_{\Omega_{\gamma}}(\alpha)$, then how the Bloch-Bohr radius will vary with respect to the domain $\Omega_{\gamma}$ for different values of $\gamma \in [0,1)$? 
\end{qsn}

\begin{qsn} \label{him-vasu-P5-qsn-1.8}
If we consider general simply connected domain $\Omega$ instead of $\Omega_{\gamma}$ which contains $\mathbb{D}$, then how the Bloch-Bohr radius for $\mathcal{B}_{\Omega}(\alpha)$ will vary with respect to the domain $\Omega$?
\end{qsn}
More precisely, it is natural to ask the following interesting questions.

\begin{qsn}  \label{him-vasu-P5-qsn-1.9}
Let $\Omega_{1} \subseteq \Omega_{2} \subseteq \cdots \subseteq \Omega_{n}$  be $n$ proper simply connected domains for some $n \geq 1$ and each domain contains $\mathbb{D}$. Let $r_{1}, r_{2}, \ldots , r_{n}$ be the Bloch-Bohr radius for $\mathcal{B}_{\Omega_{1}}(\alpha), \mathcal{B}_{\Omega_{2}}(\alpha), \ldots , \mathcal{B}_{\Omega_{n}}(\alpha)$ respectively. Then what is the relation among the Bloch-Bohr radii $r_{1}, r_{2}, \ldots , r_{n}$?  
\end{qsn}
\begin{qsn}  \label{him-vasu-P5-qsn-1.9-a}
	Let $\Omega_{1} \subseteq \Omega_{2} \subseteq \cdots \subseteq \Omega_{n}$  be $n$ proper simply connected domains for some $n \geq 1$ and each domain contains $\mathbb{D}$. Let $\tilde{r_{1}}, \tilde{r_{2}}, \ldots , \tilde{r_{n}}$ be the Bloch-Bohr radius for $\mathcal{B}_{\mathcal{H},\Omega_{1}}(\alpha), \mathcal{B}_{\mathcal{H},\Omega_{2}}(\alpha), \ldots , \mathcal{B}_{\mathcal{H},\Omega_{n}}(\alpha)$ respectively. Then what is the relation among the Bloch-Bohr radii $\tilde{r_{1}}, \tilde{r_{2}}, \ldots , \tilde{r_{n}}$?  
\end{qsn}
In this paper, we answer Question \ref{him-vasu-P5-qsn-1.9} and Question \ref{him-vasu-P5-qsn-1.9-a} completely.
\section{Affine Invariance and Inclusion Relations}
In this section, we study affine invariance property of $\mathcal{B}_{\mathcal{H},\Omega}(\alpha)$ and $\mathcal{B}^*_{\mathcal{H},\Omega}(\alpha)$. 
Let $\mathcal{L}$ be a family of harmonic mappings in $\Omega$. The family $\mathcal{L}$ is said to be affine invariant if $F(z)=a f +b \overline{f} \in \mathcal{L}$ for each $f \in \mathcal{L}$, where $a, b \in \mathbb{C}$. Further, we discuss inclusion relations for Bloch spaces on the disk $\Omega_{\gamma}$ and on arbitrary simply connected domain $\Omega$ under some suitable conditions.
\begin{prop}
The families $\mathcal{B}_{\mathcal{H},\Omega}(\alpha)$ and $\mathcal{B}^*_{\mathcal{H},\Omega}(\alpha)$ are affine invariant.
\end{prop}
\begin{pf}
Let $f \in \mathcal{B}_{\mathcal{H},\Omega}(\alpha)$ be given by $f=h + \overline{g}$. Then for $a,b \in \mathbb{C}$, we have
$$
F= af +b \overline{f}=ah+bg + \overline{\overline{a}g + \overline{b}h} := H+ \overline{G},
$$
where $H=ah+bg$ and $G=\bar{a}g + \bar{b}h$.
Therefore, \begin{align*}
|H'(z)| + |G'(z)|&= |ah'(z)+bg'(z)| + |\overline{a} g'(z)+ \overline{b}h'(z)| \\
& \leq (|a|+|b|) (|h'(z)|+|g'(z)|), \,\,z\in \Omega.
\end{align*}
It is easy to see that 
\begin{equation} \label{h-v-p5-e-2.2}
\frac{|H'(z)| + |G'(z)|}{\lambda^{\alpha}_{\Omega}(z)}\leq \frac{(|a|+|b|) (|h'(z)| + |g'(z)|)}{\lambda^{\alpha}_{\Omega}(z)},\,\,\, z\in \Omega,
\end{equation}
where $\lambda_{\varOmega}(z)$ is the hyperbolic density at $z \in \Omega$ and $\alpha \in (0,\infty)$.
Hence, from \eqref{h-v-p5-e-2.2}, $f \in \mathcal{B}_{\mathcal{H},\Omega}(\alpha)$ implies that $F= af +b \overline{f} \in \mathcal{B}_{\mathcal{H},\Omega}(\alpha)$.
On the other hand, the Jacobian of $F$ is given by $J_{F}(z)= |H'(z)|^2 - |G'(z)|^2= (|a|^2 - |b|^2)J_{f}(z)$. Thus, the conformal metrics $\lambda _{F}= \sqrt{|J_{F}|}$ and $\lambda _{f}= \sqrt{|J_{f}|}$ are homothetic {\it i.e.,} $\lambda _{F}=c \lambda _{f}$ for some $c>0$. Therefore, if $f \in \mathcal{B}^*_{\mathcal{H},\Omega}(\alpha)$ then $F= af +b \overline{f} \in \mathcal{B}^*_{\mathcal{H},\Omega}(\alpha) $.
\end{pf}

For each $\alpha>0$, both $\mathcal{B}_{\Omega_{\gamma}}(\alpha)$ and $\mathcal{B}_{\mathcal{H},\Omega_{\gamma}}(\alpha)$ are Banach spaces. In the following example, we see that $\mathcal{B}^*_{\mathcal{H},\Omega_{\gamma}}(\alpha)$ is not a linear space, which shows that some functions in $\mathcal{B}^*_{\mathcal{H},\Omega_{\gamma}}(\alpha)$ may grow arbitrarily fast. Therefore, in order to study certain properties of functions in $\mathcal{B}^*_{\mathcal{H},\Omega_{\gamma}}(\alpha)$, we shall restrict harmonic mappings to be sense-preserving. We consider the following example. The following example reduces to \cite[Example 1]{Liu-Results-Math-2018} when $\Omega_{\gamma}=\mathbb{D}$.
\begin{example}
Let $f=h+\overline{h}$, where 
$$
h(z)=\frac{\left(1-((1-\gamma)z+\gamma)\right)^{1-\beta}}{\beta -1}, \,\, z\in \Omega_{\gamma}
$$
and for some $\beta > 2 \alpha +1$. It is not difficult to see that $f$ and identity function $I(z)=z$ belong to $\mathcal{B}^*_{\mathcal{H},\Omega_{\gamma}}(\alpha)$. But we see that $F(z)=f(z)+z$ does not belong to $\mathcal{B}^*_{\mathcal{H},\Omega_{\gamma}}(\alpha)$. Indeed, 
$$
J_{F}(z)= |h'(z)+1|^2 - |h'(z)|^2 = 1+ 2 \real(h'(z)).
$$
Now we wish to find the hyperbolic density $\lambda_{\Omega}$ on $\Omega_{\gamma}$. For that we observe that $\phi: \Omega_{\gamma} \rightarrow \mathbb{D}$ define by $\phi(z)= (1-\gamma)z + \gamma$ is a conformal mapping in $\Omega_{\gamma}$. Thus, by the definition \eqref{h-v-p5-e-1.1}, we obtain the hyperbolic density for $\Omega_{\gamma}$ at the point $z$ is 
\begin{equation*}
\lambda_{\Omega_{\gamma}}(z)= \frac{1-\gamma}{1-|(1-\gamma)z + \gamma|^2}, \,\,\,\, z \in \Omega_{\gamma}.
\end{equation*}
Therefore, for $0<x<1$, we have 

\begin{align*}
\frac{|J_{F}(x)|}{\lambda^{2 \alpha} _{\Omega_{\gamma}}(x)}
&=
\frac{\left(1-((1-\gamma)x+\gamma)^2\right)^{2 \alpha}}{(1-\gamma)^{2 \alpha}} \left(1+ 2 (1-\gamma) \left(1-((1-\gamma)x+\gamma))^{-\beta}\right)\right)\\
&= \frac{\left(1+((1-\gamma)x+\gamma)\right)^{2 \alpha}}{(1-\gamma)^{2 \alpha}} \left(\frac{\left(1-((1-\gamma)x+\gamma))\right)^{\beta} +2 (1-\gamma)}{\left(1-((1-\gamma)x+\gamma))\right)^{\beta - 2 \alpha}}\right),
\end{align*}
which tends to infinity as $x \rightarrow 1 ^{-}$.
\end{example}
In the next result, we discuss inclusion relations on Bloch spaces. We shall make use of the Comparison Principle for hyperbolic metrics, which we state here only for simply connected region in complex plane. This Principle allows us to estimate the hyperbolic metric of regions in terms of other hyperbolic metrics which are known, or can be estimated easily. It is worth mentioning that it is not always possible to explicitly calculate the density of hyperbolic metric,  therefore estimates are useful.
 \begin{thm} 
 	\cite[Theorem $8.1$, Comparison Principle]{beardon-minda-hyperbolic-density}  \label{him-vasu-P5-thm-2.1-comparison-principle}
Suppose that $\Omega_{1} $ and $\Omega_{2}$ are proper simply connected regions in $\mathbb{C}$. If $\Omega_{1} \subseteq \Omega_{2}$, then $\lambda_{\Omega_{2}}(z) \leq \lambda_{\Omega_{1}}(z)$ on $\Omega_{1}$. Further, if $\lambda_{\Omega_{1}}(z) = \lambda_{\Omega_{2}}(z)$ at any point $z$ of $\Omega_{1}$, then $\Omega_{1}=\Omega_{2}$ and $\lambda_{\Omega_{1}}=\lambda_{\Omega_{2}}$.
\end{thm}
In a better way, the Comparison Principle demonstrates that the hyperbolic metric on a simply connected domain decreases as the region increases. Using this fact, we obtain the inclusion result $\mathcal{B}_{\mathcal{H}, \Omega_{2}}(\alpha) \subseteq \mathcal{B}_{\mathcal{H}, \Omega_{1}}(\alpha)$ whenever $\Omega_{1} \subseteq \Omega_{2}$.

\begin{prop} \label{him-vasu-P5-prop-2.3}
For $\alpha \in (0,\infty)$, we have 
\begin{enumerate}
	\item 
	 $\mathcal{B}_{\mathcal{H}, \Omega_{2}}(\alpha) \subseteq \mathcal{B}_{\mathcal{H}, \Omega_{1}}(\alpha)$ whenever $\Omega_{1} \subseteq \Omega_{2}$,
	\item  $\mathcal{B}_{\Omega_{\gamma}}(\alpha) \subset \mathcal{B}_{\mathcal{H},\Omega_{\gamma}}(\alpha) \subset \mathcal{B}^*_{\mathcal{H},\Omega_{\gamma}}(\alpha)$.
\end{enumerate}
\end{prop}

\begin{pf}
\begin{enumerate}
	\item Let $\Omega_{1}$ and $\Omega_{2}$ be two simply connected domains in $\mathbb{C}$ such that $\Omega_{1} \subseteq \Omega_{2}$. Let $f \in \mathcal{B}_{\mathcal{H}, \Omega_{2}}(\alpha)$. Then 
	\begin{equation} \label{him-vasu-P5-e-2.3-b}
	\frac{|h'(z)|+ |g'(z)|}{\lambda ^{\alpha}_{\Omega_{2}}(z)} < \infty \,\,\,\,\, \mbox{for} \,\,\,\,z \in \Omega_{2}.
	\end{equation}
	Since $\Omega_{1} \subseteq \Omega_{2}$, by the Comparison principle, we have $\lambda_{\Omega_{2}}(z) \leq \lambda_{\Omega_{1}}(z)$ for $z \in \Omega_{1}$. Thus, $\lambda ^{\alpha}_{\Omega_{2}}(z) \leq \lambda ^{\alpha}_{\Omega_{1}}(z)$, $z \in \Omega_{1}$ for each $\alpha \in (0,\infty)$ which leads to 
	\begin{equation} \label{him-vasu-P5-e-2.3-c}
    	\frac{|h'(z)|+ |g'(z)|}{\lambda ^{\alpha}_{\Omega_{1}}(z)} \leq	\frac{|h'(z)|+ |g'(z)|}{\lambda ^{\alpha}_{\Omega_{2}}(z)} \,\,\,\,\, \mbox{for} \,\,\,\, z \in \Omega_{1}.
	\end{equation}
	We observe that \eqref{him-vasu-P5-e-2.3-b} holds for $z \in \Omega_{2}$ and thus, \eqref{him-vasu-P5-e-2.3-b} also holds for $z \in \Omega_{1}$. This shows that 
	\begin{equation}
	\frac{|h'(z)|+ |g'(z)|}{\lambda ^{\alpha}_{\Omega_{1}}(z)} < \infty\,\,\,\,\, \mbox{for} \,\,\,\, z \in \Omega_{1},
	\end{equation} 
	which infers that $f \in \mathcal{B}_{\mathcal{H}, \Omega_{1}}(\alpha)$. Hence, $\mathcal{B}_{\mathcal{H}, \Omega_{2}}(\alpha) \subseteq \mathcal{B}_{\mathcal{H}, \Omega_{1}}(\alpha)$.
		
	\item For each $\alpha >0$, it is enough to find a function $f_{\alpha}$ such that $f_{\alpha} \in \mathcal{B}^*_{\mathcal{H},\Omega_{\gamma}}(\alpha)$ but $f_{\alpha} \not \in \mathcal{B}_{\mathcal{H},\Omega_{\gamma}}(\alpha)$. We consider the following one parameter family of functions 
	\begin{equation}\label{him-vasu-P5-e-2.3-a}
	F_{\alpha, t}(z)= H_{\alpha,t}(z)+ \overline{G_{\alpha,t}(z)} = h_{\alpha,t}(\phi(z)) + \overline{g_{\alpha,t}(\phi(z))}, \,\,\, \, t\in [0,1), \,\, z\in \Omega_{\gamma},
	\end{equation}
	where $h_{\alpha,t}$ and $g_{\alpha,t}$ are defined for $z \in \mathbb{D}$ by (see \cite{Liu-Results-Math-2018})
	\[ 
	h_{\alpha,t}(z)=
	\begin{cases}
	-\log\,(1-z) \quad &  \alpha=1/2, \\[2mm]
	\frac{(1-z)^{1/2 - \alpha} -1}{(\alpha - 1/2)} \quad &  \alpha \neq 1/2,
	\end{cases}
	\]
	and 
	\[ 
	g_{\alpha,t}(z)=
	\begin{cases}
	-\log\,(1-z) -(1-t)z \quad &  \alpha=1/2, \\[2mm]
	\frac{z}{1-z} + (1-t)log\, (1-z) \quad &  \alpha = 3/2,\\[2mm]
	\frac{(1-z)^{1/2 - \alpha} -1}{(\alpha - 1/2)} - (1-t)\frac{(1-z)^{3/2 - \alpha} -1}{(\alpha - 3/2)} \quad &  \beta \neq  \, 1/2, \, 3/2
	\end{cases}
	\]
	with $\phi : \Omega_{\gamma} \rightarrow \mathbb{D}$ defined by $\phi(z)=(1-\gamma)z + \gamma$. A simple computation shows that the dilation $w_{F_{\alpha, t}}(z)=t+(1-t)((1-\gamma)z+\gamma)$ for $z \in \Omega_{\gamma}$ and each $\alpha>0$. Then
	\begin{align*}
	\frac{\sqrt{|J_{F_{\alpha,t}}(z)|}}{\lambda^{\alpha} _{\Omega_{\gamma}}(z)}&= \frac{\left(1-|(1-\gamma)z+\gamma|^2\right)^{\alpha}}{(1-\gamma)^\alpha} \sqrt{|J_{F_{\alpha,t}}(z)|}\\
	&=  \frac{\left(1-|(1-\gamma)z+\gamma|^2\right)^{\alpha}}{(1-\gamma)^\alpha} \frac{(1-\gamma)}{\left|1- ((1-\gamma)z+\gamma)\right|^{\alpha +\frac{1}{2}}}\,\sqrt{1-\left|w_{F_{\alpha, t}}(z)\right|^2}\\
	&=
	\frac{\left(1-|\phi(z)|^2\right)^{\alpha}}{(1-\gamma)^\alpha} \frac{(1-\gamma)}{\left|1- \phi(z)\right|^{\alpha +\frac{1}{2}}}\,\sqrt{1-\left|w_{F_{\alpha, t}}(z)\right|^2}\\
	&=(1+|\phi(z)|)^\alpha \frac{(1-|\phi(z)|)^\alpha}{|1-\phi(z)|^\alpha} \sqrt{\frac{1-|\phi(z)|^2 -2t \real (\overline{\phi(z)}(1-\phi(z)))-t^2 |1-\phi(z)|^2}{|1-z|}}\\
	& \leq \frac{1}{(1-\gamma)^{\alpha -1}} 2^{\alpha + 1/2} \sqrt{1+t}, \,\, z\in \Omega_{\gamma},
	\end{align*}
	which shows that $F_{\alpha,t} \in \mathcal{B}^*_{\mathcal{H},\Omega_{\gamma}}(\alpha)$ for each $\alpha >0$. For $x \in (0,1)$, we have
	\begin{align*}
	\frac{|H^{\prime}_{\alpha,t}(x)|}{\lambda^{\alpha} _{\Omega_{\gamma}}(x)}&= \frac{\left(1-((1-\gamma)x+\gamma)^2\right)^{\alpha}}{(1-\gamma)^\alpha} \frac{(1-\gamma)}{\left(1- ((1-\gamma)x+\gamma)\right)^{\alpha +\frac{1}{2}}}\\
	&=
	\frac{\left(1+ ((1-\gamma)x+\gamma)\right)^{\alpha}}{(1-\gamma)^{\alpha -1}} \frac{1}{\sqrt{1-((1-\gamma)x+\gamma)}} \,\, \rightarrow \infty
	\end{align*}
	as $x \rightarrow 1^{-}$,
\end{enumerate}
 which shows that for each $\alpha >0$, $H^{\prime}_{\alpha,t} \not \in \mathcal{B}_{\Omega_{\gamma}}(\alpha)$ and hence, $F_{\alpha,t} \not \in \mathcal{B}_{\mathcal{H},\Omega_{\gamma}}(\alpha)$.
This completes the proof.
\end{pf}

The following result explains the structure of the set $\mathcal{B}^{*}_{\mathcal{H},\Omega_{\gamma}}(\alpha) \setminus \mathcal{B}_{\mathcal{H},\Omega_{\gamma}}(\alpha)$.

\begin{prop}
Let $f=h+\overline{g}$ be harmonic mapping in $\Omega_{\gamma}$. Then $f \in \mathcal{B}_{\mathcal{H},\Omega_{\gamma}}(\alpha)$ if, and only if, $f \in \mathcal{B}^{*}_{\mathcal{H},\Omega_{\gamma}}(\alpha)$ and either $h  \in \mathcal{B}_{\Omega_{\gamma}}(\alpha)$ or $g  \in \mathcal{B}_{\Omega_{\gamma}}(\alpha)$. Furthermore, we have 
$$
\mathcal{B}^{*}_{\mathcal{H},\Omega_{\gamma}}(\alpha) \setminus \mathcal{B}_{\mathcal{H},\Omega_{\gamma}}(\alpha)= \left\{f=h+\overline{g} \in \mathcal{B}^{*}_{\mathcal{H},\Omega_{\gamma}}(\alpha): h \not \in \mathcal{B}_{\Omega_{\gamma}}(\alpha) \,\, \mbox{and} \,\, g \not \in \mathcal{B}_{\Omega_{\gamma}}(\alpha)\right\}.
$$ 
\end{prop}

\begin{pf}
For any harmonic mapping $f = h+\overline{g}$ in $\Omega_{\gamma}$, we have
\begin{equation} \label{him-vasu-P5-e-2.5}
|h'(z)| \leq \sqrt{|J_{f}(z)|} + |g'(z)| \,\, \mbox{and} \,\, |g'(z)| \leq \sqrt{|J_{f}(z)|} + |h'(z)|,\,\, z \in \Omega_{\gamma}. 
\end{equation}
From \eqref{him-vasu-P5-e-2.5}, it is easy to see that if $f \in \mathcal{B}^{*}_{\mathcal{H},\Omega_{\gamma}}(\alpha)$ and one of $h$ and $g$ belongs to $\mathcal{B}_{\Omega_{\gamma}}(\alpha)$, then $f \in \mathcal{B}_{\mathcal{H},\Omega_{\gamma}}(\alpha)$.
\end{pf}

It is natural to ask the following question.
\begin{qsn} \label{him-P5-q-2.10}
If $f \in \mathcal{B}^{*}_{\mathcal{H},\Omega_{\gamma}}(\alpha)$ then does there exist a constant $M(\alpha)$, which depends only on $\alpha$, such that $f \in \mathcal{B}_{\mathcal{H},\Omega_{\gamma}}(M(\alpha))$?
\end{qsn}

In the following example, we point out that for a function $f \in \mathcal{B}^{*}_{\mathcal{H},\Omega_{\gamma}}(\alpha)$, whose Jacobian $J_{f}(z)=0$, does not belong to $\mathcal{B}_{\mathcal{H},\Omega_{\gamma}}(\alpha)$ for any $\alpha >0$. 

\begin{example}
Let $f=h+\overline{h}$, where 
$$
h(z)=\exp \left(\frac{1+\phi(z)}{1-\phi(z)}\right)\,\, \mbox{with} \,\, \phi(z)=(1-\gamma)z + \gamma,\, z\in \Omega_{\gamma}.
$$
Since $J_{f}(z)=|h'(z)|^2 - |h'(z)|^2=0$ for $z \in \Omega_{\gamma} $, $f \in \mathcal{B}^{*}_{\mathcal{H},\Omega_{\gamma}}(\alpha)$ for all $\alpha >0$.
\par
For $x \in (0,1)$, a simple computation shows that 
\begin{align*}
\frac{|h'(x)|}{\lambda^{\alpha} _{\Omega_{\gamma}}(x)}
&= \frac{\left(1-((1-\gamma)x +\gamma)^2\right)^{\alpha}}{(1-\gamma)^{\alpha}} \, e^{\frac{1+\phi(x)}{1-\phi (x)}} \, \frac{2(1-\gamma)}{(1-\phi(x))^2}\\
&=
\frac{\left(1-(\phi(x))^2\right)^{\alpha}}{(1-\gamma)^{\alpha -1}}\,\, \frac{2}{(1-\phi(x))^2} \,\, e^{\frac{1+\phi(x)}{1-\phi (x)}},\,\, \mbox{where}\,\, \phi(x)=(1-\gamma)x+\gamma\\
&=
\frac{2}{(1-\gamma)^{\alpha -1}} \,\, (1+\phi(x))^{2\alpha -2}\, \left(\left(\frac{1-\phi (x)}{1+ \phi (x)}\right)^{\alpha -2} \, e^{\frac{1+\phi(x)}{1-\phi (x)}}\right).
\end{align*}
Since $\phi(x) \rightarrow 1^{-}$ as $x \rightarrow 1^{-}$, we obtain 
$$
\frac{|h'(x)|}{\lambda^{\alpha} _{\Omega_{\gamma}}(x)} \rightarrow \infty \,\, \mbox{as} \,\, x\rightarrow 1^{-}.
$$
Therefore, $h \not \in \mathcal{B}_{\Omega_{\gamma}}(\alpha)$ and hence, $f \in \mathcal{B}_{\mathcal{H},\Omega_{\gamma}}(\alpha)$ for any $\alpha >0$.
\end{example}
When $\Omega_{\gamma}=\mathbb{D}$, the above example reduces to the example considered in \cite{Liu-Results-Math-2018}. From the above example, it is worth to mention that $f$ is not locally univalent in $\Omega_{\gamma}$. Therefore, to give an affirmative answer to the Question \ref{him-P5-q-2.10}, we need some additional conditions, namely locally univalent.

\begin{prop}
Let $f$ be a locally univalent harmonic mapping in $\Omega_{\gamma}$. If $f \in \mathcal{B}^{*}_{\mathcal{H},\Omega_{\gamma}}(\alpha)$, then $f \in \mathcal{B}_{\mathcal{H},\Omega_{\gamma}}(\alpha +1/2)$. Furthermore, the constant $1/2$ is sharp for each $\alpha>0$.
\end{prop}
\begin{pf}
Since $f \in \mathcal{B}_{\mathcal{H},\Omega_{\gamma}}(\alpha)$ (respectively $\mathcal{B}^{*}_{\mathcal{H},\Omega_{\gamma}}(\alpha)$) if, and only if, $\overline{f} \in \mathcal{B}_{\mathcal{H},\Omega_{\gamma}}(\alpha)$ (respectively $\mathcal{B}^{*}_{\mathcal{H},\Omega_{\gamma}}(\alpha)$), without loss of generality, we assume that $f=h+\overline{g}$ is sense-preserving in $\Omega_{\gamma}$. Then, we have
\begin{equation*}
g'(z)=\omega(z) h'(z) \,\, \mbox{and} \,\, J_{f}(z)=|h'(z)|^2 (1-|\omega(z)|^2)\,\, {\it i.e.}\,\, |h'(z)|=\sqrt{\frac{|J_{f}(z)|}{1-|\omega(z)|^2}},
\end{equation*}
where $\omega :\Omega_{\gamma} \rightarrow \mathbb{D}$ is the dilation of $f$.
\par
Now, we consider the function $\psi: \mathbb{D} \rightarrow \Omega_{\gamma}$ by $\psi(z)=(z-\gamma)/(1-\gamma)$ so that $\omega \circ \psi : \mathbb{D} \rightarrow \mathbb{D}$ is an analytic function. Then, in view of the Schwarz-Pick lemma, we have
$$
|\omega(\psi(z))| \leq \frac{|z|+|\omega(\psi(0))|}{1+ |\omega(\psi(0))| |z|},\,\, z \in \mathbb{D},
$$
{\it i.e.}
\begin{align*}
|\omega(\psi(z))| &\leq \frac{|z|+|\omega(\frac{-\gamma}{1-\gamma})|}{1+ |\omega(\frac{-\gamma}{1-\gamma})| |z|} = \frac{|z|+\omega_{0}}{1+\omega_{0}|z|},
\end{align*}
where $\omega_{0}= |\omega (\gamma/ (1-\gamma))|$.
Therefore, 
\begin{equation} \label{him-vasu-P5-e-2.8}
|\omega (z)| \leq \frac{|(1-\gamma)z+\gamma|+\omega_{0}}{1+\omega_{0} |(1-\gamma)z+\gamma|}=\frac{|\phi(z)| + \omega_{0}}{1+\omega_{0} |\phi(z)|},\, \, z \in \Omega_{\gamma},
\end{equation}
where $\phi(z)=(1-\gamma)z +\gamma$ for $z \in \Omega_{\gamma}$. Since $f \in \mathcal{B}^{*}_{\mathcal{H},\Omega_{\gamma}}(\alpha)$, we have 
\begin{equation} \label{him-vasu-P5-e-2.9}
\frac{\sqrt{|J_{f}(z)|}}{\lambda^{\alpha} _{\Omega_{\gamma}}(z)} \leq \beta^{*}_{\mathcal{H}, \Omega_{\gamma}} (\alpha)<\infty, \,\, z \in \Omega_{\gamma}.
\end{equation}
Using \eqref{him-vasu-P5-e-2.8} and \eqref{him-vasu-P5-e-2.9}, we obtain 
\begin{align*}
|h'(z)|&=\sqrt{\frac{|J_{f}(z)|}{1-|\omega(z)|^2}}\\
&\leq \beta^{*}_{\mathcal{H}, \Omega_{\gamma}} (\alpha) \,\,\lambda^{\alpha} _{\Omega_{\gamma}}(z)\,\, \frac{1}{\sqrt{1-|\omega(z)|^2}}\\
&= \frac{\beta^{*}_{\mathcal{H}, \Omega_{\gamma}} (\alpha) (1-\gamma)^{\alpha}}{(1-|(1-\gamma)z+\gamma|^2)^{\alpha}}\,\,\frac{1}{\sqrt{1-|\omega(z)|^2}}\\
&\leq \frac{\beta^{*}_{\mathcal{H}, \Omega_{\gamma}} (\alpha) (1-\gamma)^{\alpha}}{(1-|(1-\gamma)z+\gamma|^2)^{\alpha}} \,\, \left(1-\left(\frac{|\phi(z)| + \omega_{0}}{1+\omega_{0} |\phi(z)|}\right)^2\right)^{\frac{-1}{2}}\\
&=  \frac{\beta^{*}_{\mathcal{H}, \Omega_{\gamma}} (\alpha) (1-\gamma)^{\alpha}}{(1-|\phi(z)|^2)^{\alpha}} \,\, \frac{1+\omega_{0}|\phi(z)|}{\sqrt{(1-|\phi(z)|^2)(1-\omega^{2}_{0})}}\\
&\leq  \frac{\beta^{*}_{\mathcal{H}, \Omega_{\gamma}} (\alpha) (1-\gamma)^{\alpha}}{(1-|\phi(z)|^2)^{\alpha}} \,\, \frac{1+\omega_{0}}{\sqrt{(1-|\phi(z)|^2)(1-\omega^{2}_{0})}},\, \, \, \mbox{since}\, |\phi(z)|\leq 1\\
&= \frac{\beta^{*}_{\mathcal{H}, \Omega_{\gamma}} (\alpha) }{\sqrt{1-\gamma}}\,\, \sqrt{\frac{1+\omega_{0}}{1-\omega_{0}}}\,\, \frac{(1-\gamma)^{\alpha + \frac{1}{2}}}{(1-|\phi(z)|^2)^{\alpha + \frac{1}{2}}}\\
&= \frac{\beta^{*}_{\mathcal{H}, \Omega_{\gamma}} (\alpha) }{\sqrt{1-\gamma}}\,\, \sqrt{\frac{1+\omega_{0}}{1-\omega_{0}}}\,\, \lambda^{\alpha+ \frac{1}{2}} _{\Omega_{\gamma}}(z),\,\,\, \, z\in \Omega_{\gamma},
\end{align*}
which shows that $h \in \mathcal{B}_{\mathcal{H},\Omega_{\gamma}}(\alpha +1/2)$. Since $f$ is sense-preserving {\it i.e.} $|g'(z)|<|h'(z)|$ in $\Omega_{\gamma}$, $g$ also belongs to $\mathcal{B}_{\mathcal{H},\Omega_{\gamma}}(\alpha +1/2)$. Thus, $f \in \mathcal{B}_{\mathcal{H},\Omega_{\gamma}}(\alpha +1/2)$.
\vspace{3mm}

In order to show that the constant $1/2$ is sharp, we consider the function $F_{\alpha,0}$ given by \eqref{him-vasu-P5-e-2.3-a}. The function $F_{\alpha,0} \in \mathcal{B}^{*}_{\mathcal{H},\Omega_{\gamma}}(\alpha )$ and is sense-preserving in $\Omega_{\gamma}$. It is easy to see that $H_{\alpha} \in \mathcal{B}_{\Omega_{\gamma}}(\alpha +1/2)$, which implies $G_{\alpha,0} \in \mathcal{B}_{\Omega_{\gamma}}(\alpha +1/2)$. Hence, $F_{\alpha,0} \in \mathcal{B}_{\mathcal{H},\Omega_{\gamma}}(\alpha +1/2)$. A simple computation shows that $H_{\alpha} \not \in \mathcal{B}_{\Omega_{\gamma}}(p)$ for any $0<p<\alpha +1/2$. Thus, $F_{\alpha,0} \not \in \mathcal{B}_{\mathcal{H},\Omega_{\gamma}}(p)$. Therefore, the constant $\alpha +1/2$ is sharp for each $\alpha>0$.
\end{pf}

\section{Landau's theorem for harmonic Bloch mappings on the disk $\Omega_{\gamma}$}
The classical Landau Theorem for bounded analytic functions states that if $f$ is analytic in $\mathbb{D}$ with the normalizations $f(0)=0$ and $f'(0)=1$ such that $|f(z)|<M$ in $\mathbb{D}$, then $f$ is univalent in the disk $\mathbb{D}_{\rho}:=\{z: |z|<\rho\}$ with $\rho= 1/(M+ \sqrt{M^2 -1})$ and $f(\mathbb{D}_{\rho})$ contains a disk $\mathbb{D}_{R}$ with $R=M \rho^2$ (see \cite{M-S-Liu-bloch-2009}). In $1984$, Fern\'{a}ndez \cite{Fernandez-1984} studied extensively the coefficient estimate for Bloch functions. In $1987$, Colonna \cite{colona-1987} characterized the bounded analytic functions in the unit disk $\mathbb{D}$ with the Bloch functions coneections with M\"{o}bius trasformations and infinite Blaschke product. Chen {\it et al.} \cite{chen-PAMS-2000} have obtained an analogue of the Landau theorem for bounded harmonic mappings in $\mathbb{D}$. Several authors have considered Landau-type theorems for harmonic
mappings afterwards and improved their result (see \cite{chen-2014-JMMA,kalaj-2014,M-S-Liu-Landau-2009,M-S-Liu-bloch-2009,zhu-CAOT-2015,zhu-2016}). Landau's Theorem has also been extended for the classes of biharmonic mappings (see \cite{abdulhadi-20,chen-appl-math-comp-2009}). Several authors have investigated the Landau-Bloch type theorems for polyharmonic mappings (see \cite{Bai-CAOT-2019,chen-2014-JMMA}). 
We have the following inequality due to Ruscheweyh \cite{Ruscheweyh-1985}.
\begin{lem} \cite{Ruscheweyh-1985} \label{him-vasu-P5-lem-3.1}
For an analytic function $f:\mathbb{D} \rightarrow \mathbb{D}$, we have 
\begin{equation} \label{him-vasu-P5-e-3.2}
\frac{|f^{(n)}(a)|}{n!} \leq \frac{1-|f(a)|^2}{(1-|a|)^{n-1} (1-|a|^2)}
\end{equation}
for each $n\geq 1$ and $a \in \mathbb{D}$. Moreover, for each fixed $n\geq 1$ and $a \in \mathbb{D}$,
$$
\sup \limits _{f} (1-|a|)^{n-1} \,\, \frac{|f^{(n)}(a)|}{n!} \,\, \frac{1-|a|^2}{1-|f(a)|^2}=1,
$$
where the supremum is taken over all nonconstant analytic functions $f:\mathbb{D} \rightarrow \mathbb{D}$.
\end{lem}
We shall make use of the inequality \eqref{him-vasu-P5-e-3.2} to prove the following lemma.
\begin{lem} \label{him-vasu-P5-lem-3.3}
Let $f_{1}=h_{1}+ \overline{g_{1}}$ be a harmonic mapping in $\mathbb{D}$ with 
$$
h_{1}(z)=\sum \limits _{n=1}^{\infty} \alpha_{n} (z-\gamma)^n \,\, \mbox{and} \,\, g_{1}(z)=\sum \limits _{n=1}^{\infty} \beta_{n} (z-\gamma)^n \,\, \mbox{for}\, |z-\gamma|<1-\gamma .
$$
If $\lambda_{f}(\gamma)=\beta$ for some $\beta \in (0,1]$ and $||f||_{\mathcal{H},\mathbb{D}}(\alpha) \leq  M (1-\gamma)^{\alpha}$ for $M>0$, then 
$$
|\alpha_{n}|+|\beta_{n}| \leq \inf \limits  _{0<r<1} \mu (r) \,\, \mbox{for} \,\, n\geq 2,
$$
where 
\begin{equation*}
\mu(r)= \frac{m^{2}(r,\gamma)-\beta ^{2}}{nr^{n-1} (1-\gamma)^{n-1}(1-\gamma ^2) m(r, \gamma)} \,\, \mbox{and} \,\, m(r,\gamma)= \frac{M (1-\gamma)^{\alpha}}{\left(1-(1-\gamma)^2 r^2 \right)^{\alpha}}.
\end{equation*}
\end{lem}

\begin{pf}
For a fixed $r\in (0,1)$, let $F(\xi)=r^{-1} f(r(\xi - \gamma))$. Observe that $F(\xi)=H(\xi) + \overline{G(\xi)}$, where 
\begin{equation*}
H(\xi)=\frac{h_{1}(r(\xi - \gamma))}{r} \,\,\,\, \mbox{and}\,\,\,\, G(\xi)=\frac{g_{1}(r(\xi - \gamma))}{r}, \,\,\, \xi \in \mathbb{D}.
\end{equation*}
Then $F$ has the following form 
\begin{equation} \label{him-vasu-P5-e-3.4}
F(\xi) = \sum_{n=1}^{\infty} r^{n-1} \alpha_{n} (\xi - \gamma)^n + \sum_{n=1}^{\infty} r^{n-1} \overline{\beta _{n}} \,\,  \overline{(\xi - \gamma)}^n \,\, \mbox{for}\, |\xi-\gamma|<1-\gamma \, \mbox{and} \,\, \xi \in \mathbb{D}.
\end{equation}
It is easy to see that $H'(\xi)=h'(r(\xi - \gamma))$ and $G'(\xi)=g'(r(\xi - \gamma))$ for $\xi \in \mathbb{D}$.
\par
In view of the assumption $||f||_{\mathcal{H},\mathbb{D}}(\alpha) \leq M(1-\gamma)^{\alpha}$, we have 
\begin{equation*}
\frac{|h^{\prime}_{1}(z)| + |g^{\prime}_{1}(z)|}{\lambda^{\alpha}_{\mathbb{D}}(z)} \leq  M (1-\gamma)^{\alpha}
\end{equation*}
{\it i.e.},
\begin{equation}\label{him-vasu-P5-e-3.5}
|h^{\prime}_{1}(z)| + |g^{\prime}_{1}(z)| \leq M (1-\gamma)^{\alpha} \lambda^{\alpha}_{\mathbb{D}}(z)= M \, \frac{(1-\gamma)^{\alpha}}{\left(1-|z|^2\right)^{\alpha}} , \,\, z \in \mathbb{D}.
\end{equation}
 Using \eqref{him-vasu-P5-e-3.5}, we obtain
\begin{align} \label{him-vasu-P5-e-3.6}
|H'(\xi) + |G'(\xi)|
&\leq M \, \frac{(1-\gamma)^{\alpha}}{\left(1-r^{2}|\xi -\gamma|^2\right)^{\alpha}},\,\,\, \xi \in \mathbb{D},\\ \nonumber
&\leq M \, \frac{(1-\gamma)^{\alpha}}{\left(1-r^{2}(1-\gamma)^2\right)^{\alpha}}, \,\,\, \mbox{using } \,\, |\xi - \gamma|<1-\gamma,\\ \nonumber
& :=m(r, \gamma),\,\, \, \xi \in \mathbb{D}.
\end{align}
For an arbitrary $\epsilon$ with $|\epsilon|=1$, we set
$$
\widetilde{F}(\xi)=\frac{H'(\xi) +\epsilon G'(\epsilon)}{m(r,\gamma)}, \,\,\, \xi \in \mathbb{D},
$$
where $m(r,\gamma)$ is given by \eqref{him-vasu-P5-e-3.6}.
From \eqref{him-vasu-P5-e-3.4}, we see that 
\begin{equation*}
\widetilde{F}(\xi)= \frac{1}{m(r,\gamma)} \sum_{n=1}^{\infty} n (\alpha_{n} + \epsilon \beta_{n})r^{n-1} (\xi - \gamma)^{n-1}\,\, \mbox{for} \,\, |\xi - \gamma|<1-\gamma,
\end{equation*}
which is analytic in $\mathbb{D}$ and by \eqref{him-vasu-P5-e-3.6}, $|\widetilde{F}(\xi)| \leq 1$ in $\mathbb{D}$. Thus, in view of Lemma \ref{him-vasu-P5-lem-3.1}, we obtain
\begin{align} \label{him-vasu-P5-e-3.7}
\frac{nr^{n-1}|\alpha_{n} + \epsilon \beta_{n}|}{m(r,\gamma)}
&\leq \frac{1-|\widetilde{F}(\gamma)|^2}{(1-\gamma)^{n-1}(1-\gamma ^2)}\\ \nonumber
&= \frac{1- \frac{|\alpha_{1} + \epsilon \beta_{1}|^2}{m^{2}(r,\gamma)}}{(1-\gamma)^{n-1}(1-\gamma ^2)} \\ \nonumber
& \leq  \frac{1- \frac{\lambda ^{2}_{f}(\gamma)}{m^{2}(r,\gamma)}}{(1-\gamma)^{n-1}(1-\gamma ^2)} \\ \nonumber
&=  \frac{1- \frac{\beta ^{2}}{m^{2}(r,\gamma)}}{(1-\gamma)^{n-1}(1-\gamma ^2)}.
\end{align}
Since $\epsilon\, (|\epsilon|=1)$ is arbitrary, using \eqref{him-vasu-P5-e-3.7}, we deduce that 
\begin{equation*}
|\alpha_{n}| + |\beta_{n}| \leq \frac{m^{2}(r, \gamma) - \beta ^{2}}{n r^{n-1} (1-\gamma)^{n-1}(1-\gamma ^2) m(r,\gamma)}:=\mu(r)\,\, \,\, \mbox{for}\,\, n\geq 2.
\end{equation*}
For $n>1$, a simple computation shows that 
$$
\lim \limits _{r\rightarrow 0^{+}} \mu (r)= \lim \limits _{r\rightarrow 1^{-}} \mu (r)= +\infty,
$$
which ensures that infimum of $\mu (r)$ exists in $(0,1)$.
Thus, $|\alpha_{n}| + |\beta_{n}| \leq \inf \limits  _{0<r<1} \mu (r)$ for $ n\geq 2$.
\end{pf}

Using Lemma \ref{him-vasu-P5-lem-3.3}, we obtain the following coefficient estimates for $\alpha$-Bloch harmonic mappings.
\begin{thm} \label{him-vasu-P5-thm-3.1}
Let $f=h+\overline{g}$ be a harmonic mapping in $\Omega_{\gamma}$ such that 
\begin{equation*}
h(z)=\sum \limits_{n=1}^{\infty}  a_{n} z^{n} \,\,\, \mbox{and} \,\,\, g(z)=\sum \limits_{n=1}^{\infty}  b_{n} z^{n} \,\,\, \mbox{in} \,\, \mathbb{D}.
\end{equation*}
If $\lambda_{f}(0)=\lambda$ and $||f||_{\mathcal{H}, \Omega_{\gamma}} (\alpha) \leq L$, then 
\begin{equation}
|a_{n}|+|b_{n}| \leq C_{n,M}(\alpha, \lambda, \gamma)\,\,\,\, \mbox{for} \,\, n \geq 2,
\end{equation}
where 
$$
C_{n,M}(\alpha, \lambda, \gamma)=\inf \limits _{0<r<1} \frac{m^{2}(r, \gamma)-\frac{\lambda^{2}}{(1-\gamma)^2}}{n r^{n-1}(1+\gamma) m(r,\gamma)}\,\, \mbox{and}\,\, m(r,\gamma)= \frac{L}{1-\gamma}\,\,  \frac{ (1-\gamma)^{\alpha}}{\left(1-(1-\gamma)^2 r^2 \right)^{\alpha}}.
$$
\end{thm}

\begin{pf}
Let $\psi : \mathbb{D} \rightarrow \Omega_{\gamma}$ be defined by  $\psi(z)=(z-\gamma)/(1-\gamma)$. Then the composition $T\circ \psi$ is harmonic in $\mathbb{D}$ such that $T(z)=h(\psi(z))+\overline{g(\psi(z))}$ with 
\begin{equation}
T(z)= \sum_{n=1}^{\infty} \frac{a_{n}}{(1-\gamma)^{n}} (z-\gamma)^n + \overline{\sum_{n=1}^{\infty} \frac{b_{n}}{(1-\gamma)^{n}} (z-\gamma)^n}  \,\, \, \mbox{for} \,\, |z-\gamma|<1-\gamma ,\,\,\, z \in \mathbb{D}.
\end{equation}
Set 
$$
h_{T}(z)=h\left(\frac{z-\gamma}{1-\gamma}\right) \,\, \mbox{and} \,\, g_{T}(z)=g\left(\frac{z-\gamma}{1-\gamma}\right).
$$
Then $T(z)= h_{T}(z) + \overline{g_{T}(z)}$ for $z \in \mathbb{D}$. Clearly,
\begin{equation} \label{him-vasu-P5-e-3.10}
h^{\prime}_{T}(z)= \frac{1}{1-\gamma} h^{\prime}\left(\frac{z- \gamma}{1-\gamma}\right) \,\, \mbox{and} \,\, g^{\prime}_{T}(z)= \frac{1}{1-\gamma} g^{\prime}\left(\frac{z- \gamma}{1-\gamma}\right), \,\,\, z \in \mathbb{D}.
\end{equation}
Since $\lambda_{f}(0)=\lambda$, we obtain 
\begin{equation} \label{him-vasu-P5-e-3.11}
\lambda_{T}(\gamma)= \left||h^{\prime}_{T}(\gamma)|- |g^{\prime}_{T}(\gamma)|\right|=\frac{1}{1-\gamma} ||h'(0)|-|g'(0)||=\frac{\lambda_{f}(0)}{1-\gamma}=\frac{\lambda}{1-\gamma}=\beta \,\,(\mbox{say}).
\end{equation}
In view of the assumption $||f||_{\mathcal{H}, \Omega_{\gamma}} (\alpha) \leq L$, we have 
\begin{equation} \label{him-vasu-P5-e-3.12}
|h'(z)|+|g'(z)| \leq \frac{L(1-\gamma)^{\alpha}}{\left(1-|(1-\gamma)z+\gamma|^2\right)^{\alpha}}, \,\, \, z \in \Omega_{\gamma}.
\end{equation}
Then, from \eqref{him-vasu-P5-e-3.10} and \eqref{him-vasu-P5-e-3.12}, we obtain 
\begin{equation} \label{him-vasu-P5-e-3.13}
|h^{\prime}_{T}(z)|+ |g^{\prime}_{T}(z)| \leq \frac{1}{1-\gamma} \frac{L(1-\gamma)^{\alpha}}{(1-|z|^{2})^{\alpha}}, \,\, z \in \mathbb{D}.
\end{equation}
In view of \eqref{him-vasu-P5-e-3.11}, \eqref{him-vasu-P5-e-3.13} and setting $\beta = \lambda/(1-\gamma)$ and $M=L/(1-\gamma)$, applying Lemma \ref{him-vasu-P5-lem-3.3} to the function $T$, we obtain
 \begin{align*}
 \frac{1}{(1-\gamma)^n} (|a_{n}|+|b_{n}|) 
 &\leq \inf \limits _{0<r<1} \frac{m^{2}(r, \gamma)-\frac{\lambda^{2}}{(1-\gamma)^2}}{n r^{n-1}(1-\gamma)^{n-1}(1-\gamma ^{2}) m(r,\gamma)},
 \end{align*}
 where 
 $$
 m(r,\gamma)= \frac{L}{1-\gamma}\,\,  \frac{ (1-\gamma)^{\alpha}}{\left(1-(1-\gamma)^2 r^2 \right)^{\alpha}} , \,\, r\in (0,1).
 $$
 Therefore,
 \begin{align*}
 |a_{n}|+|b_{n}| &\leq \inf \limits _{0<r<1} \frac{m^{2}(r, \gamma)-\frac{\lambda^{2}}{(1-\gamma)^2}}{n r^{n-1}(1+\gamma) m(r,\gamma)} := C_{n,M}(\alpha,\gamma,\lambda).
 \end{align*}
 This completes the proof.
\end{pf}

By making use of Theorem \ref{him-vasu-P5-thm-3.1}, we establish the Landau's Theorem for the bounded normalized functions in $\mathcal{B}_{\mathcal{H},\Omega _{\gamma}}(\alpha)$.
\begin{thm}
Let $f=h+\overline{g}$ be a harmonic mapping in $\Omega_{\gamma}$ with $f(0)=\lambda_{f}(0)-\lambda=0$ and $||f||_{\mathcal{H}, \Omega_{\gamma}} (\alpha)\leq M$ such that 
\begin{equation*}
h(z)=\sum \limits_{n=1}^{\infty}  a_{n} z^{n} \,\,\, \mbox{and} \,\,\, g(z)=\sum \limits_{n=1}^{\infty}  b_{n} z^{n} \,\,\, \mbox{in} \,\, \mathbb{D}.
\end{equation*}
Then $f$ is univalent in $\mathbb{D}_{\rho_{0}}:=\{z:|z|<\rho_{0}\}$, where $\rho_{0}$ is the unique root of 
\begin{equation}
\sum_{n=2}^{\infty} C_{n,M}(\alpha, \lambda, \gamma) \, n r^{n-1}=\lambda
\end{equation}
in $(0,1)$. Moreover, $f(\mathbb{D}_{\rho_{0}})$ contains the disk $\mathbb{D} _{\rho}$, where 
$$
\rho = \lambda \, \rho_{0} - \sum_{n=2}^{\infty} C_{n,M}(\alpha,\lambda,\gamma) \rho_{0} ^n.
$$
\end{thm}

\begin{pf}
Let $z_{1}, z_{2} \in \mathbb{D}_{r}$, $r<1$, where $0<r<r_{\rho_{0}}$ and $z_{1} \not = z _{2}$. Since $\lambda_{f}(0)=\lambda$, we have
\begin{equation} \label{him-vasu-P5-e-3.15}
\lambda_{f}(0)=||a_{1}|-|b_{1}||=\lambda .
\end{equation}
In view of Theorem \ref{him-vasu-P5-thm-3.1}, we have 
\begin{equation} \label{him-vasu-P5-e-3.16}
|a_{n}| + |b_{n}| \leq C_{n,M}(\alpha, \lambda, \gamma) \,\,\, \mbox{for} \,\, n \geq 2.
\end{equation}
By making use of \eqref{him-vasu-P5-e-3.15} and \eqref{him-vasu-P5-e-3.16}, we obtain 
\begin{align} \label{him-vasu-P5-e-3.17}
|f(z_{1}) - f(z_{2})| &= \left|\,\, \int \limits _{[z_{1},z_{2}]} f_{z}(z)\, dz + f_{\overline{z}}(z)\, d \overline{z}\right|\\[2mm]  \nonumber 
&\geq 
\left|\,\, \int \limits _{[z_{1},z_{2}]} f_{z}(0)\, dz + f_{\overline{z}}(0)\, d \overline{z}\right|\\[2mm] \nonumber
& - \left|\,\, \int \limits _{[z_{1},z_{2}]} (f_{z}(z)-f_{z}(0))\, dz +  (f_{\overline{z}}(z) - f_{\overline{z}}(0))\, d \overline{z}\right|\\[2mm] \nonumber
&> |z_{1}-z_{2}|\, \lambda_{f}(0)- |z_{1}-z_{2}| \sum_{n=2}^{\infty} (|a_{n}|+|b_{n}|)n r^{n-1} \\[2mm] \nonumber
&= |z_{1}-z_{2}| \left[\lambda - \sum_{n=2}^{\infty} C_{n,M}(\alpha, \lambda, \gamma) n r^{n-1}\right]\\[2mm] \nonumber
&:= |z_{1} - z_{2}|\,  \Psi(r).
\end{align}
Now see that $\Psi$ is continuous and differentiable in $(0,1)$ and 
$$\Psi ^{'}(r)=-\sum_{n=2}^{\infty} C_{n,M}(\alpha, \lambda, \gamma)\, n\, (n-1)\, r^{n-2} < 0
$$
 for $r \in (0,1)$, which shows that $\Psi$ is decreasing in $(0,1)$. It is important to note that $\Psi (0)=\lambda >0$ and $\lim _{r \rightarrow 1 ^{-}} \Psi(r)= - \infty$, which ensures that $\Psi$ has the unique root in $(0,1)$ and choose that to be $\rho_{0} $. Then, from \eqref{him-vasu-P5-e-3.17}, we deduce that $|f(z_{1})-f(z_{2})|>0$ for any $0<r<\rho_{0}$. This shows that $f$ is univalent in $\mathbb{D}_{\rho_{0}}$.
\vspace{3mm}

Clearly, $f(0)=0$. Then, for $z'=\rho_{0} e^{i\theta} \in \partial \mathbb{D}_{\rho_{0}}$, we have 
\begin{align} \label{him-vasu-P5-e-3.18}
|f(z')| &\geq |a_{1}z' + \overline{b_{1}}\,\,\overline{z'}|- \left|\sum_{n=2}^{\infty} (a_{n}z'^{n} + \overline{b_{n}}\,\,\overline{z'}^n)\right|\\ \nonumber
&\geq \lambda_{f}(0)\rho_{0} - \sum_{n=2}^{\infty} (|a_{n}|+|b_{n}|) \rho_{0}^{n} \\ \nonumber
& \geq \lambda \, \rho_{0} - \sum_{n=2}^{\infty} C_{n,M}(\alpha,\lambda,\gamma) \rho_{0} ^n.
\end{align}
Therefore, \eqref{him-vasu-P5-e-3.18} shows that $f(\mathbb{D}_{\rho_{0}})$ contains the disk $\mathbb{D} _{\rho}$, where 
$$
\rho = \lambda \, \rho_{0} - \sum_{n=2}^{\infty} C_{n,M}(\alpha,\lambda,\gamma) \rho_{0} ^n.
$$
This completes the proof.
\end{pf}

In the next Theorem, we obtain the coefficient estimates for the bounded functions in $\mathcal{B}_{\mathcal{H},\Omega}(\alpha)$. In particular, we obtain the coefficient estimates when $\Omega=\Omega_{\gamma}$ and $\Omega=\mathbb{D}$.
\begin{thm}
Let $f=h+\overline{g}$ be a harmonic mapping in $\Omega$ such that $f(0)=0$ and $||f||_{\mathcal{H}, \Omega} (\alpha)\leq M$ for some constant $M>0$, where 
\begin{equation*}
h(z)=\sum \limits_{n=1}^{\infty}  a_{n} z^{n} \,\,\, \mbox{and} \,\,\, g(z)=\sum \limits_{n=1}^{\infty}  b_{n} z^{n} \,\,\, \mbox{in} \,\, \mathbb{D}.
\end{equation*}
Then the following inequality 
\begin{equation}
|a_{n}|^2 + |b_{n}|^2 \leq  A_{n}(\Omega, M)
\end{equation}
holds for all $n \geq 1$, where 
$$
A_{n}(\Omega, M) := \frac{M^2}{n^2}\,\, \inf _{0<t<1} \frac{1}{2\pi \,t\, t^{2(n-1)}} \int \limits _{|z|=t} \lambda ^{2 \alpha}_{\Omega}(z) \,\, |d\,z|.
$$
Moreover, $|a_{n}| + |b_{n}| \leq \sqrt{2 A_{n}(\Omega, M)}$ for $n \geq 1$. In particular, we have
\begin{enumerate}
	\item If $\Omega=\Omega_{\gamma}$, then $|a_{n}|^2 + |b_{n}|^2 \leq  C_{n}(\alpha, \gamma, M)$ for $n \geq 1$, where 
	\begin{equation*}
	C_{n}(\alpha, \gamma, M):= \frac{M^2}{n^2}\,\, \inf _{0<t<1} \frac{(1-\gamma)^{2 \alpha}}{t^{2(n-1)} (1-((1-\gamma)t + \gamma)^2)^{2 \alpha}}. 
	\end{equation*}

	\item  If $\Omega= \mathbb{D}$, then $|a_{n}|^2 + |b_{n}|^2 \leq  C_{n}(\alpha, 0, M):=C_{n}(\alpha, M)$ for $n \geq 1$ and  \\ $\lim  _{n \rightarrow \infty} C_{n}(\alpha, 0, M)=0$ for each $\alpha \in (0, 1)$. If $\alpha \geq 1$, then 
	\begin{equation*}
	C_{n}(\alpha, M) \leq \frac{M^2}{(2\alpha)^{2\alpha}}\,\, \left(1+ \frac{2\alpha}{n-1}\right)^{n-1}\,\, \frac{(n-1+2\alpha)^{2\alpha}}{n^2}\,= {\it O}(n^{2\alpha - 2}).
	\end{equation*}
\end{enumerate}
\end{thm}

\begin{pf} 
The assumptions $f(0)=0$ and $||f||_{\mathcal{H}, \Omega} (\alpha)\leq M$	 show that 
\begin{equation*}
|h'(z)| + |g'(z)| \leq M \lambda^{\alpha}_{\Omega}(z), \,\,\, z \in \Omega.
\end{equation*}
It is easy to see that 
\begin{equation} \label{him-vasu-P5-e-3.20}
|h'(z)|^2 + |g'(z)|^2 \leq (|h'(z)| + |g'(z)|)^2 \leq M^2\, \, \lambda^{2\alpha}_{\Omega}(z), \,\, \, z \in \Omega.
\end{equation}
Integrating the inequality \eqref{him-vasu-P5-e-3.20} over the circle $|z|=t<1$, we obtain 
\begin{equation} \label{him-vasu-P5-e-3.21}
2\, \pi t \sum_{n=1}^{\infty} n^2 (|a_{n}|^2 + |b_{n}|^2)t^{2(n-1)} \leq  M^2\, \, \int \limits _{|z|=t} \lambda^{2\alpha}_{\Omega}(z)\, |dz|, \,\, \, z \in \mathbb{D}.
\end{equation}
Thus, from \eqref{him-vasu-P5-e-3.21}, it follows that 
$$
|a_{n}|^2 + |b_{n}|^2 \leq \frac{M^2}{n^2}\,\, \frac{1}{2\pi \,t\, t^{2(n-1)}} \int \limits _{|z|=t} \lambda ^{2 \alpha}_{\Omega}(z) \,\, |d\,z| =\frac{M^2}{n^2}\,\, \Psi _{1}(t),
$$
where 
$$
\Psi _{1}(t)= \frac{1}{2\pi \,t\, t^{2(n-1)}} \int \limits _{|z|=t} \lambda ^{2 \alpha}_{\Omega}(z) \,\, |d\,z|.
$$
Thus, we now only need to show that $\inf_{0<t<1} \Psi_{1}(t)$ exists. To prove this, we make use of the comparison principle of hyperbolic density function to the domains $\mathbb{D}$ and $\Omega$: if $\mathbb{D} \subseteq \Omega$ then 
\begin{equation} \label{him-vasu-P5-e-3.22}
\lambda_{\Omega}(z) \leq \lambda_{\mathbb{D}}(z), \,\,\, z \in \mathbb{D}.
\end{equation}
The inequality \eqref{him-vasu-P5-e-3.22} leads to 
\begin{equation} \label{him-vasu-P5-e-3.23}
\Psi_{1}(t)=\frac{1}{2\pi \,t\, t^{2(n-1)}} \int \limits _{|z|=t} \lambda ^{2 \alpha}_{\Omega}(z) \,\, |d\,z| \leq \frac{1}{2\pi \,t\, t^{2(n-1)}} \int \limits _{|z|=t} \lambda ^{2 \alpha}_{\mathbb{D}}(z) \,\, |d\,z|:= \Psi_{2}(t), \,\,\, z \in \mathbb{D}.
\end{equation}
A simple computation using the fact $\lambda _{\mathbb{D}}(z)=1/(1-|z|^2)$, shows that 
\begin{align} \label{him-vasu-P5-e-3.24}
\Psi_{2}(t) &= \frac{1}{2\pi \,t\, t^{2(n-1)}} \int \limits _{|z|=t}  \frac{|d\,z|}{(1-|z|^2)^{2 \alpha}}\\ \nonumber
&= \frac{1}{2\pi \,t\, t^{2(n-1)}} \int \limits _{\theta=0}^ {2 \pi} \frac{t}{(1-|te^{i \theta}|^2)^{2 \alpha}}\,\, d \theta \\ \nonumber
&= \frac{1}{2\pi \,\, t^{2(n-1)}} \frac{2 \pi}{(1-t^2)^{2 \alpha}} \\ \nonumber
&= \frac{1}{t^{2(n-1)}(1-t^2)^{2 \alpha}}.
\end{align}
From \eqref{him-vasu-P5-e-3.24}, we note that 
$$
\lim\limits_{t \rightarrow 0} \Psi_{2} (t)=\lim\limits_{t \rightarrow 1^{-}} \Psi_{2} (t)= +\infty,
$$
which ensures that $\inf_{0<t<1} \Psi_{2}(t)$ exists and hence, by \eqref{him-vasu-P5-e-3.23}, $\inf_{0<t<1} \Psi_{1}(t)$ exists.
Therefore, we have
\begin{equation*}
|a_{n}|^2 + |b_{n}|^2 \leq  A_{n}(\Omega, M) \,\,n \geq 1,
\end{equation*}
where 
\begin{equation} \label{him-vasu-P5-e-3.25}
A_{n}(\Omega, M) = \frac{M^2}{n^2}\,\, \inf _{0<t<1} \frac{1}{2\pi \,t\, t^{2(n-1)}} \int \limits _{|z|=t} \lambda ^{2 \alpha}_{\Omega}(z) \,\, |d\,z|.
\end{equation}
\vspace{2mm}

We note that 
\begin{equation*}
|a_{n}| + |b_{n}| \leq \sqrt{2(|a_{n}|^2 + |b_{n}|^2)} \leq \sqrt{2 A_{n}(\Omega, M)} \,\,\,\, \mbox{for} \,\, n \geq 1.
\end{equation*}

If $\Omega = \Omega_{\gamma}$, then 
\begin{equation} \label{him-vasu-P5-e-3.26}
\lambda_{\Omega}(z) = \lambda_{\Omega_{\gamma}}(z)= \frac{(1-\gamma)^{2 \alpha}}{(1-|(1-\gamma)z + \gamma|^2)^{2 \alpha}} \leq \frac{(1-\gamma)^{2 \alpha}}{(1-((1-\gamma)|z| + \gamma)^2)^{2 \alpha}}. 
\end{equation}

Using \eqref{him-vasu-P5-e-3.21} and \eqref{him-vasu-P5-e-3.26}, we obtain 

\begin{align} \label{him-vasu-P5-e-3.27}
2\, \pi t \sum_{n=1}^{\infty} n^2 (|a_{n}|^2 + |b_{n}|^2)t^{2(n-1)} 
&\leq M^2 \int \limits _{|z|=t} \frac{(1-\gamma)^{2 \alpha}}{(1-((1-\gamma)|z| + \gamma)^2)^{2 \alpha}}\,\, |dz| \\ \nonumber
 &=  M^2 \,\, \int \limits _{\theta=0}^{2\pi}\frac{ t (1-\gamma)^{2 \alpha}}{(1-((1-\gamma)|te^{i\theta}| + \gamma)^2)^{2 \alpha}}\,\, d\theta \\ \nonumber
 &= M^2 \frac{2\pi t (1-\gamma)^{2 \alpha}}{(1-((1-\gamma)t + \gamma)^2)^{2 \alpha}}
\end{align}
and therefore, we can express \eqref{him-vasu-P5-e-3.27} as 
\begin{equation} \label{him-vasu-P5-e-3.28}
|a_{1}|^2 + |b_{1}|^2 + 2^2 (|a_{1}|^2 + |b_{1}|^2)t^2 + \cdots \leq \frac{M^2 (1-\gamma)^{2\alpha}}{(1-\gamma ^2)^{2 \alpha}} + A_{1}t + \cdots.
\end{equation}
Thus, from \eqref{him-vasu-P5-e-3.28}, we deduce that 
$$
|a_{1}|^2 + |b_{1}|^2 \leq \frac{M^2}{(1+\gamma)^{2\alpha}}.
$$

Using \eqref{him-vasu-P5-e-3.25}, we have 
\begin{align} \label{him-vasu-P5-e-3.29}
A_{n}(\Omega_{\gamma}, M)&= \frac{M^2}{n^2}\,\, \inf _{0<t<1} \frac{1}{2\pi \,t\, t^{2(n-1)}} \int \limits _{|z|=t} \frac{(1-\gamma)^{2 \alpha}}{(1-|(1-\gamma)z + \gamma|^2)^{2 \alpha}}\,\, |dz| \\ \nonumber
& \leq \frac{M^2}{n^2}\,\, \inf _{0<t<1} \frac{1}{2\pi \,t\, t^{2(n-1)}} \int \limits _{|z|=t} \frac{(1-\gamma)^{2 \alpha}}{(1-((1-\gamma)|z| + \gamma)^2)^{2 \alpha}}\,\, |dz| \\ \nonumber
& = \frac{M^2}{n^2}\,\, \inf _{0<t<1} \frac{(1-\gamma)^{2 \alpha}}{t^{2(n-1)} (1-((1-\gamma)t + \gamma)^2)^{2 \alpha}} \\ \nonumber
& =:C_{n}(\alpha, \gamma, M).
\end{align}
Let 
$$
\mu(t)= \frac{(1-\gamma)^{2 \alpha}}{t^{2(n-1)} (1-((1-\gamma)t + \gamma)^2)^{2 \alpha}},\,\, t \in (0,1).
$$
For $n \geq 2$, we can see that
$$
\lim\limits_{t \rightarrow 0} \mu (t)= \lim\limits_{t \rightarrow 1^{-}} \mu (t)=+ \infty.
$$
This observation shows that the infimum of $\mu(t)$ must exists in $(0,1)$. For $n \geq 2$, we compute that 
$$
\mu '(t)= (1-\gamma)^{2 \alpha} \left(\frac{1}{t^{2(n-1)}} \, \frac{4 \alpha (1-\gamma) ((1-\gamma)t + \gamma)}{(1-((1-\gamma)t + \gamma)^2)^{2 \alpha + 1}} - \frac{2(n-1)}{t^{2(n-1) +1}} \,\, \frac{1}{(1-((1-\gamma)t + \gamma)^2)^{2 \alpha }}\right).
$$
For each $0<\alpha <1$,  $\mu ' (t)=0$ has the following roots
$$
t_{1}= \frac{-\gamma(n-1+\alpha)+ \sqrt{\alpha ^2 \gamma ^2 + (n-1) (n-1 + 2\alpha)}}{(1-\gamma)(n-1+2\alpha)} <1
$$
and 
$$
t_{2}= \frac{-\gamma(n-1+\alpha)- \sqrt{\alpha ^2 \gamma ^2 + (n-1) (n-1 + 2\alpha)}}{(1-\gamma)(n-1+2\alpha)}<0
$$
such that $t_{1} \in (0,1)$ and $t_{2}<0$ for each $\alpha \in (0, \infty)$, and thus, we obtain that
\begin{equation*}
\inf _{0<t<1} \mu (t)= \mu (t_{1}).
\end{equation*}
It follows from \eqref{him-vasu-P5-e-3.29} that 
\begin{equation} \label{him-vasu-P5-e-3.30}
C_{n}(\alpha, \gamma, M) \leq \frac{M^2}{n^2} \mu (t_{1}).
\end{equation} 
\vspace{3mm}

We observe that $\Omega_{\gamma}$ reduces to $\mathbb{D}$ for $\gamma=0$. Then $t_{1}= \sqrt{(n-1)/(n-1+2\alpha)}$ and 
\begin{align} \label{him-vasu-P5-e-3.31}
C_{n}(\alpha,0, M)&:=C_{n}(\alpha,M) \\ \nonumber
& = \frac{M^2}{n^2} \left(1+ \frac{2\alpha}{n-1}\right)^{n-1}\,\, \frac{(n-1+2\alpha)^{2\alpha}}{(2\alpha)^{2\alpha}} \\ \nonumber
&= \frac{M^2}{(2\alpha)^{2\alpha}} \,\, \left(1+ \frac{2\alpha}{n-1}\right)^{n-1}\,\, \frac{(n-1+2\alpha)^{2\alpha}}{n^2}.
\end{align}
We see that if $0<\alpha<1$, $C_{n}(\alpha,M) \rightarrow 0$ as $n \rightarrow + \infty$. If $\alpha=1$, then 
$$
C_{n}(\alpha,M) = \frac{M^2}{4}\,\, \left(1+ \frac{2}{n-1}\right)^{n-1}\,\, \left(\frac{n+1}{n}\right)^2.
$$
If $\alpha >1$, then 
$$
C_{n}(\alpha,M) \leq \frac{M^2}{(2\alpha)^{2\alpha}} \,\, \left(1+ \frac{2\alpha}{n-1}\right)^{n-1}\,\, \frac{(n-1+2\alpha)^{2\alpha}}{n^2}= {\it O}\left(n^{2\alpha - 2}\right).
$$
This completes the proof.
\end{pf}

\section{Bohr inequality for Bloch spaces}
In this section, we study Bloch-Bohr radius for the Bloch spaces. We first obtain the Bloch-Bohr radius for $\mathcal{B}_{\Omega}(\alpha)$, when $\Omega$ is arbitrary proper simply connected domain in $\mathbb{C}$.
\begin{thm} \label{him-vasu-P5-thm-4.1}
Let $\Omega$ be a proper simply connected domain containing $\mathbb{D}$. Let $f \in \mathcal{B}_{\Omega}(\alpha)$ with $||f||_{\Omega, \alpha} \leq 1$ such that $f(z)= \sum _{n=0}^{\infty} a_{n} z^{n}$ in $\mathbb{D}$. Then $\sum _{n=0}^{\infty}  |a_{n}| r^n \leq 1$
for $|z|=r \leq r_{\Omega}(\alpha)$, where $r_{\Omega}(\alpha)$ is the smallest root of 
\begin{equation*}
I(r):=\frac{r}{2 \pi}\, \int \limits _{|z|=r} \lambda ^{2 \alpha}_{\Omega}(z)\,\, |dz| = \frac{6}{\pi ^{2}}
\end{equation*}
in $(0,1)$, provided $\lim _{r \rightarrow 1^{-}}I(r)>6/\pi^2$.
\end{thm}

\begin{pf}
Let $f \in \mathcal{B}_{\Omega}(\alpha)$ with $||f||_{\Omega, \alpha} \leq 1$. Then 
\begin{equation*}
|f(0)| + \sup \limits_{z \in \Omega} \frac{|f'(z)|}{\lambda ^{\alpha}_{\Omega}(z)} \leq 1,
\end{equation*}	
which implies that $|f'(z)| \leq (1-|a_{0}|)\, \lambda ^{\alpha}_{\Omega}(z) $, $z \in \Omega$. Therefore, 
\begin{equation} \label{him-vasu-P5-e-4.1}
|f'(z)|^2 \leq (1- |a_{0}|)^2 \, \lambda ^{2\alpha}_{\Omega}(z)\,\,\,\, \mbox{for} \,\,\,\, z \in \Omega.
\end{equation}
Since $f(z)= \sum _{n=0}^{\infty} a_{n} z^{n}$ in $\mathbb{D}$, from \eqref{him-vasu-P5-e-4.1}, we obtain 
\begin{equation} \label{him-vasu-P5-e-4.2}
\left|\sum \limits_{n=1}^{\infty} n \, a_{n} \, z^{n-1}\right|^2 \leq 
(1- |a_{0}|)^2 \, \lambda ^{2\alpha}_{\Omega}(z)\,\,\,\, \mbox{for} \,\,\,\, z \in \mathbb{D}.
\end{equation}
Integrating \eqref{him-vasu-P5-e-4.2} over the circle $|z|=r<1$, we obtain 
\begin{equation*}
2 \pi r \sum \limits_{n=1}^{\infty} n^2 \, |a_{n}|^2 \, r^{2(n-1)} \leq (1-|a_{0}|)^2 \int \limits _{|z|=r} \lambda ^{2 \alpha}_{\Omega}(z)\,\, |dz|\,\,\,\, \mbox{for} \,\,\,\, z \in \mathbb{D},
\end{equation*}
which leads to 
\begin{equation} \label{him-vasu-P5-e-4.3}
\sum \limits_{n=1}^{\infty} n^2 |a_{n}|^2 r^{2n} \leq
(1-|a_{0}|)^2 \, \frac{r}{2 \pi}\, \int \limits _{|z|=r} \lambda ^{2 \alpha}_{\Omega}(z)\,\, |dz| \,\,\,\,\mbox{for} \,\, z \in \mathbb{D}.
\end{equation}
In view of the classical Cauchy-Schwartz inequality and \eqref{him-vasu-P5-e-4.3}, we obtain
\begin{align} \label{him-vasu-P5-e-4.4}
|a_{0}| + \sum \limits_{n=1}^{\infty} |a_{n}| r^{n} & \leq 
|a_{0}| + \sqrt{\sum \limits_{n=1}^{\infty} n^2 |a_{n}|^2 r^{2n}}\,\, \sqrt{\frac{\pi ^2}{6}} \\ \nonumber
& \leq |a_{0}| + (1-|a_0|) \sqrt{\frac{r}{2 \pi}\, \int \limits _{|z|=r} \lambda ^{2 \alpha}_{\Omega}(z)\,\, |dz|} \,\, \sqrt{\frac{\pi ^2}{6}}.
\end{align}
We note that the right hand side of \eqref{him-vasu-P5-e-4.4} less than or equals to $1$ if 
\begin{equation*}
\sqrt{\frac{r}{2 \pi}\, \int \limits _{|z|=r} \lambda ^{2 \alpha}_{\Omega}(z)\,\, |dz|} \,\, \sqrt{\frac{\pi ^2}{6}} \leq 1, 
\end{equation*}
 or equivalently, if
\begin{equation} \label{him-vasu-P5-e-4.5}
\frac{r}{2 \pi}\, \int \limits _{|z|=r} \lambda ^{2 \alpha}_{\Omega}(z)\,\, |dz| \leq \frac{6}{\pi ^{2}},
\end{equation}
which holds for $r \leq r_{\Omega}(\alpha)$, where $r_{\Omega}(\alpha)$ is the smallest root of 
\begin{equation*}
I(r):=\frac{r}{2 \pi}\, \int \limits _{|z|=r} \lambda ^{2 \alpha}_{\Omega}(z)\,\, |dz| = \frac{6}{\pi ^{2}}
\end{equation*}
in $(0,1)$. In order to prove the existence of the root $r_{\Omega}(\alpha)$ in $(0,1)$, we consider the function $H: [0,1) \rightarrow \mathbb{R}$ defined by 
\begin{equation*}
H(r)= I(r)- \frac{6}{\pi ^{2}} = \frac{r}{2 \pi}\, \int \limits _{|z|=r} \lambda ^{2 \alpha}_{\Omega}(z)\,\, |dz| - \frac{6}{\pi ^{2}}.
\end{equation*}
It is easy to see that $H$ is continuous in $[0,1)$ with 
\begin{equation} \label{him-vasu-P5-e-4.6}
H(0)=I(0)- \frac{6}{\pi ^{2}}=-\frac{6}{\pi ^{2}}<0
\end{equation}
and
\begin{equation} \label{him-vasu-P5-e-4.7}
\lim \limits _{r \rightarrow 1 ^{-}} H(r)= \lim \limits _{r \rightarrow 1 ^{-}} I(r) - \frac{6}{\pi ^{2}}>0 \,\,\, \, (\mbox{by the assumption}\,\, \lim \limits _{r \rightarrow 1 ^{-}} I(r) > 6/\pi ^{2}).
\end{equation}
Therefore, in view of \eqref{him-vasu-P5-e-4.6} and \eqref{him-vasu-P5-e-4.7}, we conclude that $H$ has a root in $(0,1)$ and choose the smallest root to be $r_{\Omega}(\alpha)$. Thus, from \eqref{him-vasu-P5-e-4.4}, we deduce that $\sum_{n=0}^{\infty} |a_{n}|r^n \leq 1$ for $r \leq r_{\Omega}(\alpha)$. This completes the proof.
\end{pf}

In the next result, we independently obtain the Bloch-Bohr radius for $\mathcal{B}_{\Omega_{\gamma}}(\alpha)$.
\begin{thm} \label{him-vasu-P5-thm-4.2}
For $0 \leq \gamma <1$, let $f \in \mathcal{B}_{\Omega_{\gamma}}(\alpha)$ with $||f||_{\Omega_{\gamma}, \alpha} \leq 1$ such that $f(z)= \sum _{n=0}^{\infty} a_{n} z^{n}$ in $\mathbb{D}$. Then $\sum _{n=0}^{\infty}  |a_{n}| r^n \leq 1$
for $|z|=r \leq r_{\gamma}(\alpha)$, where $r_{\gamma}(\alpha)$ is the unique root of $H_{\gamma, \alpha}(r)=0$ in $(0,1)$, where
\begin{equation*}
H_{\gamma, \alpha}(r):=6(1-((1-\gamma)r + \gamma)^2)^{2\alpha} - \pi ^2 \, (1-\gamma)^{2 \alpha}\, r^2.
\end{equation*}
Moreover, the radius $r_{\gamma}(\alpha)$ cannot be replaced by a number greater than $R_{\gamma}(\alpha)$ when $\alpha>1$, where 
\begin{equation*}
R_{\gamma}(\alpha)= \frac{1}{(1-\gamma)}\,\, \left(-\gamma + \sqrt{1- (1-\gamma)((1+\gamma)^{1-\alpha} + 2(\alpha -1))^{\frac{1}{1-\alpha}}}\,\,\right).
\end{equation*}
\end{thm}

\begin{pf}
The hyperbolic density for $\Omega_{\gamma}$ at $z$ is 
\begin{equation*}
\lambda_{\Omega_{\gamma}}(z)= \frac{1-\gamma}{1-|(1-\gamma)z + \gamma|^2}, \,\,\,\, z \in \Omega_{\gamma}.
\end{equation*}
Then, by the hypothesis $||f||_{\Omega_{\gamma}, \alpha} \leq 1$, we obtain
\begin{align*}
|f'(z)| &\leq (1-|a_{0}|) \,  \frac{(1-\gamma)^{\alpha}}{(1-|(1-\gamma)z + \gamma|^2)^{\alpha}}\\
& \leq (1-|a_{0}|) \, \frac{(1-\gamma)^{\alpha}}{(1-((1-\gamma)|z| + \gamma)^2)^{\alpha}}, \,\,\, z \in \Omega_{\gamma}
\end{align*}
and thus, 
\begin{equation} \label{him-vasu-P5-e-4.8}
\left|\sum \limits_{n=1}^{\infty} n \, a_{n} \, z^{n-1}\right|^2 =|f'(z)|^2 \leq  (1-|a_{0}|)^2 \, \frac{(1-\gamma)^{2\alpha}}{(1-((1-\gamma)|z| + \gamma)^2)^{2\alpha}}, \,\,\,\, z \in \mathbb{D}.
\end{equation}
Integrating \eqref{him-vasu-P5-e-4.8} over the circle $|z|=r<1$, we obtain 
\begin{align*}
2 \pi r \sum \limits_{n=1}^{\infty} n^2 \, |a_{n}|^2 \, r^{2(n-1)} & \leq (1-|a_{0}|)^2 \,\, \int \limits _{|z|=r} \frac{(1-\gamma)^{2\alpha}}{(1-((1-\gamma)|z| + \gamma)^2)^{2\alpha}} \, |dz| \\
& = (1-|a_{0}|)^2 \,\, \int \limits_{\theta=0}^{2\pi} \frac{(1-\gamma)^{2\alpha}}{(1-((1-\gamma)|re^{i \theta}| + \gamma)^2)^{2\alpha}} \, r d \theta \\
&= (1-|a_{0}|)^2 \,\, \frac{(1-\gamma)^{2\alpha}}{(1-((1-\gamma)r + \gamma)^2)^{2\alpha}}\, 2\pi r,
\end{align*} 
which is equivalent to 
\begin{equation} \label{him-vasu-P5-e-4.9}
\sum \limits_{n=1}^{\infty} n^2 \, |a_{n}|^2 \, r^{2n} \leq 
(1-|a_{0}|)^2 \,\, \frac{(1-\gamma)^{2\alpha}\, r^2}{(1-((1-\gamma)r + \gamma)^2)^{2\alpha}}.
\end{equation}
Using Cauchy-Schwartz inequality and \eqref{him-vasu-P5-e-4.9}, we obtain
\begin{align} \label{him-vasu-P5-e-4.10}
|a_{0}| + \sum \limits_{n=1}^{\infty} |a_{n}| r^{n} & \leq 
|a_{0}| + \sqrt{\sum \limits_{n=1}^{\infty} n^2 |a_{n}|^2 r^{2n}}\,\, \sqrt{\frac{\pi ^2}{6}} \\ \nonumber
& \leq |a_{0}| + (1-|a_{0}|)\,  \frac{(1-\gamma)^{\alpha}\, r}{(1-((1-\gamma)r + \gamma)^2)^{\alpha}} \,\, \sqrt{\frac{\pi ^2}{6}} \leq 1,
\end{align}
provided 
\begin{equation} \label{him-vasu-P5-e-4.11}
\frac{(1-\gamma)^{\alpha}\, r}{(1-((1-\gamma)r + \gamma)^2)^{\alpha}} \,\, \sqrt{\frac{\pi ^2}{6}} \leq 1.
\end{equation}
Now \eqref{him-vasu-P5-e-4.11} holds for $r \leq r _{\gamma}(\alpha)$, where $r_{\gamma}(\alpha) \in (0,1)$ is the smallest root of 
\begin{equation*}
\frac{(1-\gamma)^{\alpha}\, r}{(1-((1-\gamma)r + \gamma)^2)^{\alpha}} \,\, \sqrt{\frac{\pi ^2}{6}} = 1,
\end{equation*}
or equivalently 
\begin{equation*}
6(1-((1-\gamma)r + \gamma)^2)^{2\alpha} - \pi ^2 \, (1-\gamma)^{2 \alpha}\, r^2=0.
\end{equation*}
To prove the existence of the root $r_{\gamma}(\alpha)$, we consider the function $H_{\gamma, \alpha}:[0,1] \rightarrow \mathbb{R}$ defined by 
\begin{equation*}
H_{\gamma, \alpha}(r)=6(1-((1-\gamma)r + \gamma)^2)^{2\alpha} - \pi ^2 \, (1-\gamma)^{2 \alpha}\, r^2.
\end{equation*} 
Clearly, $H_{\gamma, \alpha}$ is continuous in $[0,1]$ and satisfies the conditions
\begin{equation} \label{him-vasu-P5-e-4.12}
H_{\gamma, \alpha}(0)= 6 (1-\gamma ^2)^{2\alpha}>0 \,\, \, \mbox{and} \,\,\, H_{\gamma, \alpha}(1)= -\pi ^2 \, (1-\gamma)^{2 \alpha}\, r^2<0.
\end{equation}
Thus, applying the Intermediate value theorem to the continuous function $H_{\gamma, \alpha}$ in $[0,1]$, we can easily see that $H_{\gamma, \alpha}$ has a root in $(0,1)$ and choose the smallest root to be $r_{\gamma}(\alpha)$. Hence, from  \eqref{him-vasu-P5-e-4.10}, we obtain $\sum_{n=0}^{\infty} |a_{n}| r^{n} \leq 1$ for $r \leq r_{\gamma}(\alpha)$.
\par
In order to obtain the upper bound of $r_{\gamma}(\alpha)$ when $\alpha >1$, we consider the function
\begin{equation*}
f_{\gamma, \alpha}(z)=(1-\gamma)^{\alpha -1} \,\,  \frac{(1-((1-\gamma)z+\gamma)^2)^{1-\alpha} - (1-\gamma ^2)^{1- \alpha}}{2(\alpha -1)}, \,\,\, z \in \Omega_{\gamma}.
\end{equation*}
We see that 
\begin{equation*}
f'_{\gamma,\alpha}(z)= \frac{(1-\gamma)^{\alpha -1} (\alpha - 1) (1-\gamma)((1-\gamma)z + \gamma)}{2(\alpha -1) (1-((1-\gamma)z +\gamma)^2)^{\alpha}}, \,\,\,\, z \in \Omega_{\gamma},
\end{equation*}
and hence we compute that $||f_{\gamma,\alpha}||_{\Omega_{\gamma},\alpha}=1$ when $\alpha > 1$. It is easy to see that $f_{\gamma, \alpha}= f_{\alpha} \circ \phi $ for $z \in \Omega_{\gamma}$, where $f_{\alpha}: \mathbb{D} \rightarrow \mathbb{C}$ is analytic function in $\mathbb{D}$ and defined by
\begin{equation} \label{him-vasu-P5-e-4.12-a}
f_{\alpha}(z)=\frac{(1-z^2)^{1-\alpha} -1}{2(\alpha -1)}= \sum \limits_{n=1}^{\infty} b_{\alpha,n} \, z^n, \,\,\,\,\,\, z \in \mathbb{D}
\end{equation}
and $\phi: \Omega_{\gamma} \rightarrow \mathbb{D}$ defined by $\phi(z)=(1-\gamma)z + \gamma$. We point out that all the Taylor's coefficients $b_{\alpha,n}$ of $f_{\alpha}(z)$ in \eqref{him-vasu-P5-e-4.12-a} are non-negative real numbers for $\alpha>1$ and also, the Taylor's coefficients of $\phi(z)$ are non-negative for each $0 \leq \gamma <1$. Then if $f_{\gamma,\alpha}(z)= \sum _{n=1}^{\infty} b_{\alpha,\gamma,n} \, z^n$ in $\mathbb{D}$, we observe that all the Taylor's coefficients $b_{\alpha,\gamma,n}$ of $f_{\gamma,\alpha}(z)$ are also non-negative for each $\alpha>1$ and $0 \leq \gamma <1$. Therefore, for $\alpha >1$, we obtain 
\begin{equation} \label{him-vasu-P5-e-4.13}
\sum \limits_{n=1}^{\infty} b_{\alpha,\gamma,n} \, r^n = 
(1-\gamma)^{\alpha -1} \,\,  \frac{(1-((1-\gamma)r+\gamma)^2)^{1-\alpha} - (1-\gamma ^2)^{1- \alpha}}{2(\alpha -1)}, \,\,\, \, 0<r<1.
\end{equation}
Now we wish to find the smallest $r$ such that $\sum \limits_{n=1}^{\infty} b_{\alpha,\gamma,n} \, r^n =1$,
which gives
\begin{equation} \label{him-vasu-P5-e-4.14}
(1-\gamma)^{\alpha -1} \,\,  \frac{(1-((1-\gamma)r+\gamma)^2)^{1-\alpha} - (1-\gamma ^2)^{1- \alpha}}{2(\alpha -1)}=1.
\end{equation}
A simple computation using \eqref{him-vasu-P5-e-4.14} shows that
\begin{equation*}
r= \frac{1}{(1-\gamma)}\,\, \left(-\gamma + \sqrt{1- (1-\gamma)[(1+\gamma)^{1-\alpha} + 2(\alpha -1)]^{\frac{1}{1-\alpha}}}\,\,\right):=R_{\gamma}(\alpha).
\end{equation*}
Therefore, the radius $r_{\gamma}(\alpha)$ cannot be replaced by a number greater than $R_{\gamma}(\alpha)$. 
\end{pf}

\begin{table}[ht]
	\centering
	\begin{tabular}{|l|l|l|}
		\hline
		$\alpha$& $r_{0.1}(\alpha)$& $r_{0.4}(\alpha)$ \\
		\hline
		$(0,0.5]$& $(0.779697 \searrow 0.619322]$& $(0.779697 \searrow 0.631373]$\\
		\hline
		$(0.5,1]$& $(0.619322 \searrow 0.554985]$& $(0.631373 \searrow 0.576500]$\\
		\hline
		$(1,1.5]$& $(0.554985 \searrow 0.514933]$& $(0.576500 \searrow 0.544243]$\\
		\hline
		$(1.5,2]$& $(0.514933 \searrow 0.48638]$& $(0.544243 \searrow 0.522310]$\\
		\hline
		$(2,2.5]$& $(0.48638 \searrow 0.464523]$& $(0.522310 \searrow 0.506191]$\\
		\hline
		$(2.5,3]$& $(0.464523 \searrow 0.447025]$& $(0.506191 \searrow 0.493744]$\\
		\hline
	\end{tabular}
	\vspace{3mm}
	\caption{Values of $r_{0.1}(\alpha)$ and $r_{0.4}(\alpha)$ for various values of $\alpha$.}
	\label{tabel-4.2-a}
\end{table}

\begin{figure}[!htb]
	\begin{center}
		\includegraphics[width=0.48\linewidth]{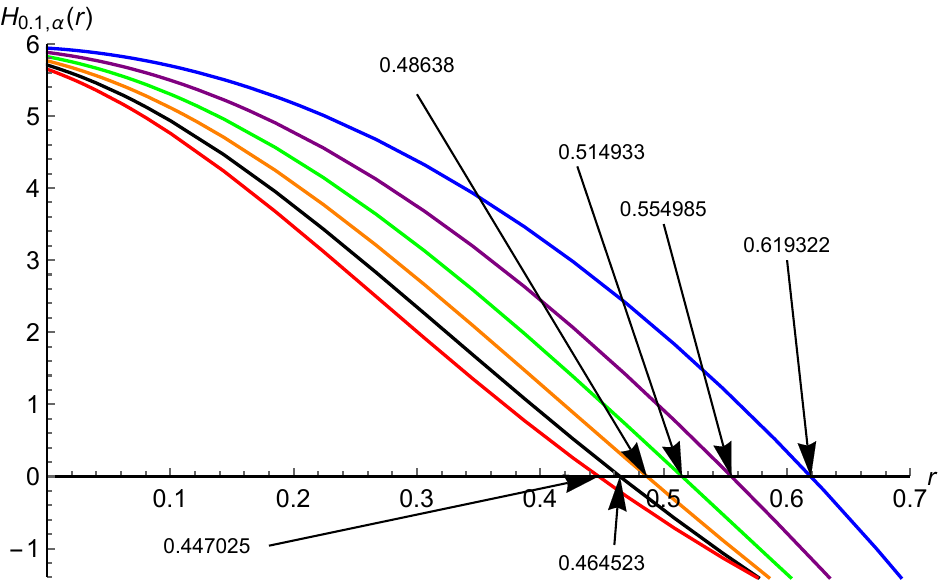}
		\,
		\includegraphics[width=0.48\linewidth]{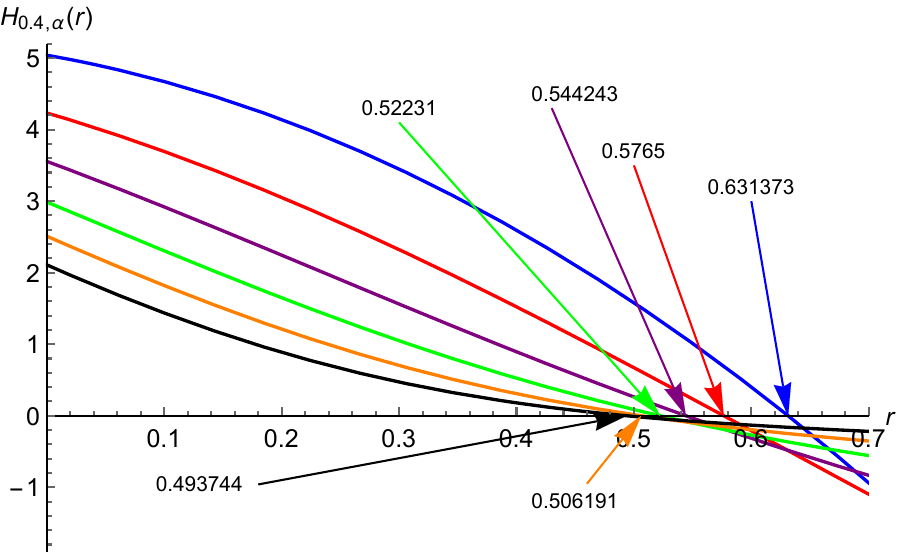}
	\end{center}
	\caption{The graph of $H_{0.1}(\alpha, r) $ and $H_{0.4}(\alpha, r)$ in $ (0,1) $ when $\alpha=0.1,0.2,0.3,0.4, 0.5, 0.6$.}
	\label{figure-4.2-a}
\end{figure}

\begin{table}[ht]
	\centering
	\begin{tabular}{|l|l|l|}
		\hline
		$\alpha$& $r_{0.7}(\alpha)$& $r_{0.9}(\alpha)$ \\
		\hline
		$(0,0.5]$& $(0.779697 \searrow 0.641889]$& $(0.779697 \searrow 0.648220]$\\
		\hline
		$(0.5,1]$& $(0.641889 \searrow 0.594235]$& $(0.648220 \searrow 0.604518]$\\
		\hline
		$(1,1.5]$& $(0.594235 \searrow 0.567402]$& $(0.604518 \searrow 0.580481]$\\
		\hline
		$(1.5,2]$& $(0.567402 \searrow 0.549746]$& $(0.580481 \searrow 0.564932]$\\
		\hline
		$(2,2.5]$& $(0.549746 \searrow 0.537110]$& $(0.564932 \searrow 0.553950]$\\
		\hline
		$(2.5,3]$& $(0.537110 \searrow 0.527566]$& $(0.553950 \searrow 0.545743]$\\
		\hline
	\end{tabular}
	\vspace{3mm}
	\caption{The Values of $r_{0.7}(\alpha)$ and $r_{0.9}(\alpha)$ for various values of $\alpha$.}
	\label{tabel-4.2-b}
\end{table}

\begin{figure}[!htb]
	\begin{center}
		\includegraphics[width=0.48\linewidth]{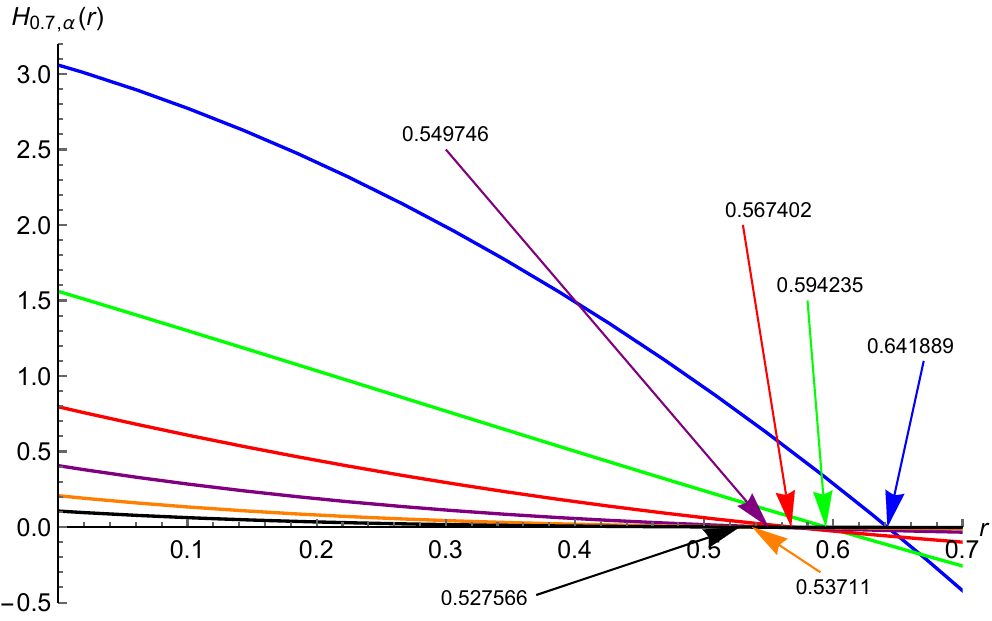}
		\,
		\includegraphics[width=0.49\linewidth]{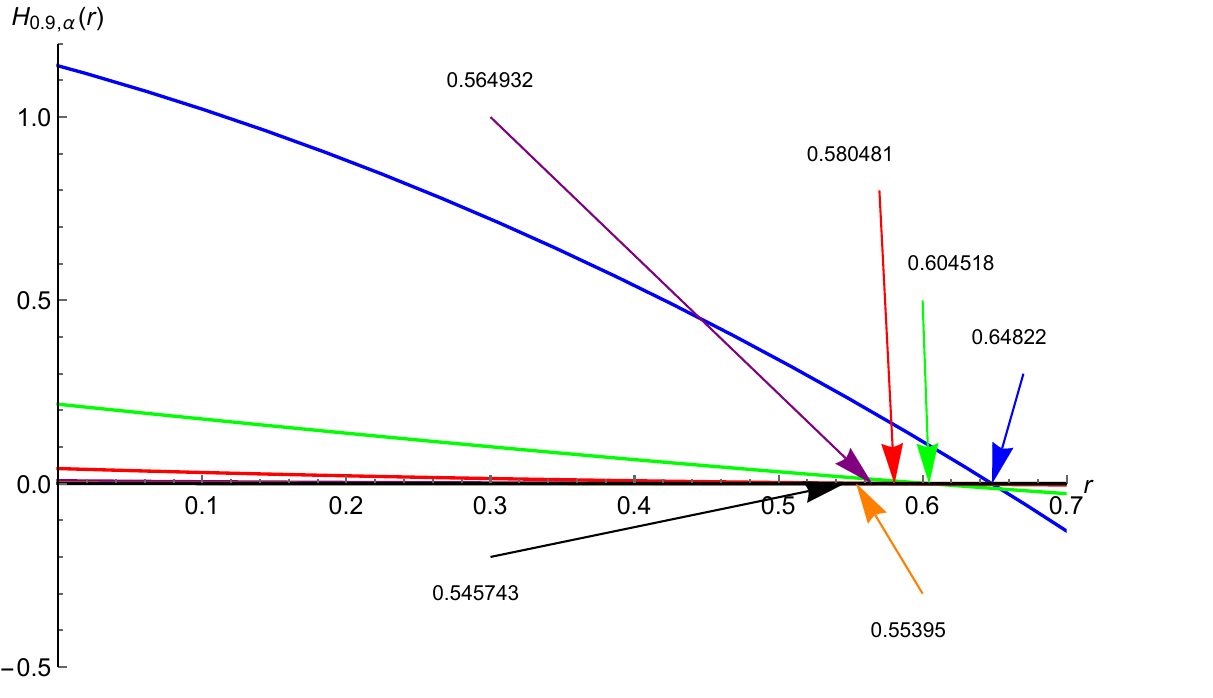}
	\end{center}
	\caption{The graph of $H_{0.7}(\alpha, r) $ and $H_{0.9}(\alpha, r)$ in $ (0,1) $ when $\alpha=0.1,0.2,0.3,0.4, 0.5, 0.6$.}
	\label{figure-4.2-b}
\end{figure}

\begin{table}[ht]
	\centering
	\begin{tabular}{|l|l|l|l|l|}
		\hline
		$\alpha$& $R_{0.1}(\alpha)$& $R_{0.4}(\alpha)$&  $R_{0.7}(\alpha)$ & $R_{0.9}(\alpha)$ \\
		\hline
		$1.5$& $0.860175$& $0.846027$& $0.835811$& $0.830626$\\
		\hline
		$2.0$& $0.812265$& $0.804300$& $0.800870$& $0.800085$\\
		\hline
		$2.5$& $0.774074$& $0.771636$& $0.773354$& $0.775558$\\
		\hline
		$3.0$& $0.742592$& $0.745067$& $0.750847$& $0.755236$\\
		\hline
	\end{tabular}
	\vspace{3mm}
	\caption{The Values of $R_{0.1}(\alpha), R_{0.4}(\alpha), R_{0.7}(\alpha)$ and $R_{0.9}(\alpha)$ for $\alpha=1.5, 2.0, 2.5, 3.0$.}
	\label{tabel-4.2-c}
\end{table}
From Table \ref{tabel-sec-1}, Table \ref{tabel-4.2-a} and Table \ref{tabel-4.2-b}, we observe that Bloch-Bohr radius is greater than Bohr radius for each $\alpha \in (0,\infty)$ and $\gamma \in [0,1)$. In the Table \ref{tabel-4.2-a} and Table \ref{tabel-4.2-b}, the notation $(r_{0.1}(\alpha_{1}) \searrow r_{0.1}(\alpha_{2})]$ means that the value of $r_{0.1}(\alpha)$ is monotonically decreasing from $\lim _{\alpha \rightarrow \alpha ^{+}_{1}}= r_{0.1}(\alpha_{1})$  to $ r_{0.1}(\alpha_{2})$ when $\alpha_{1} < \alpha \leq \alpha_{2}$.
From Table \ref{tabel-4.2-a} and Table \ref{tabel-4.2-b}, it is easy to see that the Bloch-Bohr radius is monotonically decreasing in $\alpha$. We also observe that, from Table \ref{tabel-4.2-a} and Table \ref{tabel-4.2-b}, the Bloch-Bohr radius is monotonically increasing in $\gamma$. This observation leads us to compare the Bloch-Bohr radius for $\mathcal{B}_{\Omega_{1}}(\alpha)$ and $\mathcal{B}_{\Omega_{2}}(\alpha)$, where $\Omega_{1}$ and $\Omega_{2}$ are proper simply connected domains containing $\mathbb{D}$ such that $\Omega_{1} \subseteq \Omega_{2}$.

\begin{cor} \label{him-vasu-P5-cor-4.15}
Let $\Omega_{1}$ and $\Omega_{2}$ be two proper simply connected domains in $\mathbb{C}$ containing $\mathbb{D}$ such that $\Omega_{1} \subseteq \Omega_{2}$. Suppose that Bloch-Bohr radius exists for both $\mathcal{B}_{\Omega_{1}}(\alpha)$ and $\mathcal{B}_{\Omega_{2}}(\alpha)$, and call them as $r_{1}$ and $r_{2}$ respectively. Then $r_{1} \leq r_{2}$.
\end{cor}

\begin{pf}
In view of the comparison principle, we have $\lambda ^{\alpha}_{\Omega_{2}}(z) \leq \lambda ^{\alpha}_{\Omega_{1}}(z)$ for $z \in \Omega_{1}$. Then, in view of the inequality \eqref{him-vasu-P5-e-4.4} in Theorem \ref{him-vasu-P5-thm-4.1}, we obtain 
\begin{align*} 
\sum \limits_{n=0}^{\infty} |a_{n}| r^n &\leq  |a_{0}| + (1-|a_0|) \sqrt{\frac{r}{2 \pi}\, \int \limits _{|z|=r} \lambda ^{2 \alpha}_{\Omega_{2}}(z)\,\, |dz|} \,\, \sqrt{\frac{\pi ^2}{6}} \\    \nonumber 
& \leq |a_{0}| + (1-|a_0|) \sqrt{\frac{r}{2 \pi}\, \int \limits _{|z|=r} \lambda ^{2 \alpha}_{\Omega_{1}}(z)\,\, |dz|} \,\, \sqrt{\frac{\pi ^2}{6}}
\end{align*}
for each $r \in (0,1)$. This shows that Bloch-Bohr radius $r_{1}$ for $\mathcal{B}_{\Omega_{1}}(\alpha)$ is less than or equals to Bloch-Bohr radius $r_{2}$ for $\mathcal{B}_{\Omega_{2}}(\alpha)$.
\end{pf}

As a consequence of Corollary \ref{him-vasu-P5-cor-4.15}, for simply connected domains $\Omega_{1}, \Omega_{2}, \ldots , \Omega_{n}$ such that $\Omega_{1} \subseteq \ldots \subseteq\Omega_{n}$, we obtain the following result.

\begin{cor} \label{him-vasu-P5-cor-4.15-a}
Let $\Omega_{1}, \Omega_{2}, \ldots , \Omega_{n}$ be finite $n$- proper simply connected domains in $\mathbb{C}$ containing $\mathbb{D}$ such that $\Omega_{1} \subseteq \Omega_{2}\subseteq \ldots \subseteq \Omega_{n}$. Suppose that Bloch-Bohr radius exists for all $\mathcal{B}_{\Omega_{1}}(\alpha), \mathcal{B}_{\Omega_{2}}(\alpha), \ldots, $  $\mathcal{B}_{\Omega_{n}}(\alpha)$, and call them as $r_{1}, r_{2}, \ldots , r_{n}$ respectively. Then $r_{1} \leq r_{2} \leq \ldots \leq r_{n}$.
\end{cor} 

In $2018$, Kayumov and Ponnusamy \cite{Kayumov-Ponnusamy-2018-b} extensively studied the $p$-Bohr radius for harmonic functions in $\mathbb{D}$. Motivated by $p$-Bohr radius, in this article, we study the $p$-Bloch-Bohr radius for Bloch spaces. In the next result, for $p \geq 1$, we obtain $p$-Bloch-Bohr radius for $\mathcal{B}_{\mathcal{H}, \Omega_{\gamma}}(\alpha)$. For functions $f=h+ \overline{g}$ of the form \eqref{him-vasu-P5-4.18-d}, for each $p \geq 1$, we define $p$-Bloch-Bohr radius for $\mathcal{B}_{\mathcal{H}, \Omega}(\alpha)$ to be the largest radius $r(p) \in (0,1)$ such that 
\begin{equation*}
|a_{0}| + \sum \limits_{n=1}^{\infty} \left(|a_{n}|^p + |b_{n}|^p\right)^{1/p}\, r^n \leq 1 \,\,\,\,\, \mbox{for} \,\,\, |z|=r \leq r(p)
\end{equation*} 
for all $f \in \mathcal{B}_{\mathcal{H}, \Omega}(\alpha)$. 
When $f$ is analytic we can obtain analytic $p$-Bloch-Bohr radius. In particular, for $p=1$, analytic $p$-Bloch-Bohr radius coincides with Bohr radius.
\begin{thm}
Let $f=h+\overline{g} \in \mathcal{B}_{\mathcal{H}, \Omega_{\gamma}}(\alpha)$ such that $||f||_{\mathcal{H}, \Omega_{\gamma}} (\alpha)\leq 1$ with
\begin{equation} \label{him-vasu-P5-4.18-d}
h(z)=\sum \limits_{n=1}^{\infty}  a_{n} z^{n} \,\,\, \mbox{and} \,\,\, g(z)=\sum \limits_{n=1}^{\infty}  b_{n} z^{n} \,\,\, \mbox{in} \,\, \mathbb{D}.
\end{equation}
Then, for each $p \geq 1$, we have 
\begin{equation*}
|a_{0}| + \sum \limits_{n=1}^{\infty} \left(|a_{n}|^p + |b_{n}|^p\right)^{1/p}\, r^n \leq 1
\end{equation*}
for $|z|=r \leq r_{\gamma}(\alpha,p)$, where $r_{\gamma}(\alpha,p)$ is the smallest root of $H_{1}(r)=0$ in $(0,1)$, where 
\begin{equation*}
H_{1}(r)=6(1-((1-\gamma)r + \gamma)^2)^{2\alpha} - K_{p}\,\pi ^2 \, (1-\gamma)^{2 \alpha}\, r^2
\end{equation*}
 and $K_{p}=\max \{2^{(2/p)-1}, 1\}$.
\end{thm}

\begin{pf}
Since $||f||_{\mathcal{H}, \Omega_{\gamma}} (\alpha)\leq 1$, we obtain
\begin{equation*}
|h'(z)| + |g'(z)| \leq \frac{(1-|a_0|)(1-\gamma)^{\alpha}}{(1-|(1-\gamma)z + \gamma|^2)^{\alpha}} \leq  \, \frac{(1-|a_{0}|)(1-\gamma)^{\alpha}}{(1-((1-\gamma)|z| + \gamma)^2)^{\alpha}}, \,\,\, z \in \Omega_{\gamma}
\end{equation*}
and hence  
\begin{equation} \label{him-vasu-P5-e-4.15}
|h'(z)|^2 + |g'(z)|^2 \leq (|h'(z)| + |g'(z)|)^2 \leq  \frac{(1-|a_{0}|)^2\,(1-\gamma)^{2\alpha}}{(1-((1-\gamma)|z| + \gamma)^2)^{2\alpha}}, \,\, \, z \in \Omega.
\end{equation}
Integrating the inequality \eqref{him-vasu-P5-e-4.15} over the circle $|z|=r<1$, we obtain 
\begin{equation} \label{him-vasu-P5-e-4.16}
\sum \limits_{n=1}^{\infty} n^2 (|a_{n}|^2 + |b_{n}|^2)r^{2(n-1)} \leq   \frac{(1-|a_{0}|)^2(1-\gamma)^{2\alpha}}{(1-((1-\gamma)r + \gamma)^2)^{2\alpha}}.
\end{equation}
Using the classical Cauchy-Schwartz inequality and \eqref{him-vasu-P5-e-4.16}, we obtain 
\begin{align} \label{him-vasu-P5-e-4.17}
|a_{0}| + \sum \limits_{n=1}^{\infty}\left(|a_{n}|^p + |b_{n}|^p\right)^{1/p}\, r^n & \leq 
|a_{0}| + \sqrt{\sum \limits_{n=1}^{\infty} n^2 (|a_{n}|^p + |b_{n}|^p) ^{2/p} \, r^{2n}} \, \sqrt{\sum \limits_{n=1}^{\infty} \frac{1}{n^2}} \\ \nonumber
& \leq |a_{0}| + \sqrt{K_{p} \, \sum \limits_{n=1}^{\infty} n^2 (|a_{n}|^2 + |b_{n}|^p2) \, r^{2n}}\, \sqrt{\frac{\pi ^2}{6}} \\ \nonumber
& \leq |a_{0}| + \sqrt{K_{p}} \, \frac{(1-|a_{0}|)(1-\gamma)^{\alpha}\, r}{(1-((1-\gamma)r + \gamma)^2)^{\alpha}}\, \sqrt{\frac{\pi ^2}{6}} ,
\end{align}
which is less than or equals to $1$ for $r \leq r_{\gamma}(\alpha,p)$, where $r_{\gamma}(\alpha,p)$ is the smallest root of 
\begin{equation*}
 \frac{\sqrt{K_{p}}(1-\gamma)^{\alpha}\, r}{(1-((1-\gamma)r + \gamma)^2)^{\alpha}}\, \sqrt{\frac{\pi ^2}{6}}=1,
\end{equation*}
or equivalently,
\begin{equation} \label{him-vasu-P5-e-4.18}
6(1-((1-\gamma)r + \gamma)^2)^{2\alpha} - K_{p}\,\pi ^2 \, (1-\gamma)^{2 \alpha}\, r^2=0
\end{equation}
in $(0,1)$. Using the same lines of argument as in the proof of Theorem \ref{him-vasu-P5-thm-4.2}, we can show that \eqref{him-vasu-P5-e-4.18} has a root in $(0,1)$ and choose $r_{\gamma}(\alpha,p)$ to be the smallest root in $(0,1)$. This completes the proof.
\end{pf}

\begin{table}[ht]
	\centering
	\begin{tabular}{|l|l|l|}
		\hline
		$\alpha$& $r_{0.1}(\alpha,1)$& $r_{0.4}(\alpha,1)$ \\
		\hline
		$(0,0.5]$& $(0.551329 \searrow 0.489073]$& $(0.551329 \searrow 0.505842]$\\
		\hline
		$(0.5,1]$& $(0.489073 \searrow 0.454081]$& $(0.505842 \searrow 0.482236]$\\
		\hline
		$(1,1.5]$& $(0.454081 \searrow 0.430126]$& $(0.482236 \searrow 0.467118]$\\
		\hline
		$(1.5,2]$& $(0.430126 \searrow 0.412188]$& $(0.467118 \searrow 0.456423]$\\
		\hline
		$(2,2.5]$& $(0.412188 \searrow 0.398028]$& $(0.456423 \searrow 0.448389]$\\
		\hline
		$(2.5,3]$& $(0.398028 \searrow 0.386449]$& $(0.448389 \searrow 0.442102]$\\
		\hline
	\end{tabular}
	\vspace{3mm}
	\caption{Values of $r_{0.1}(\alpha,1)$ and $r_{0.4}(\alpha,1)$ for various values of $\alpha$.}
	\label{tabel-4.3-a}
\end{table}

\begin{figure}[!htb]
	\begin{center}
		\includegraphics[width=0.48\linewidth]{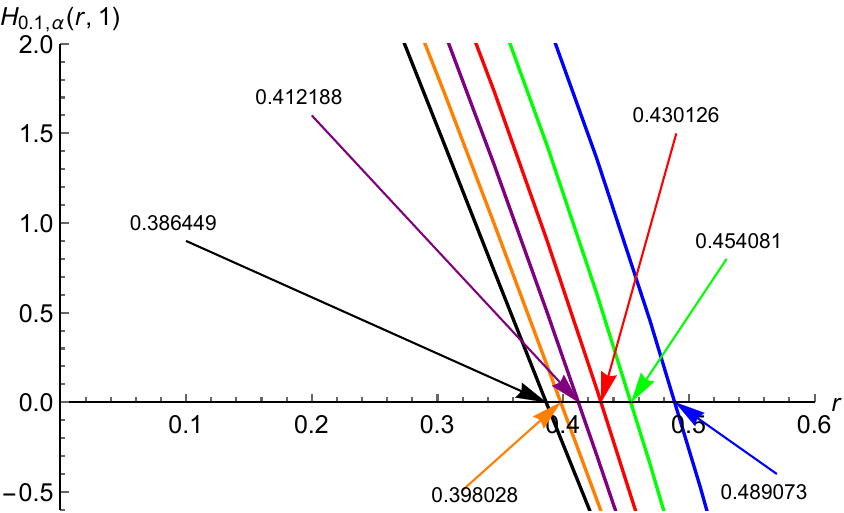}
		\,
		\includegraphics[width=0.48\linewidth]{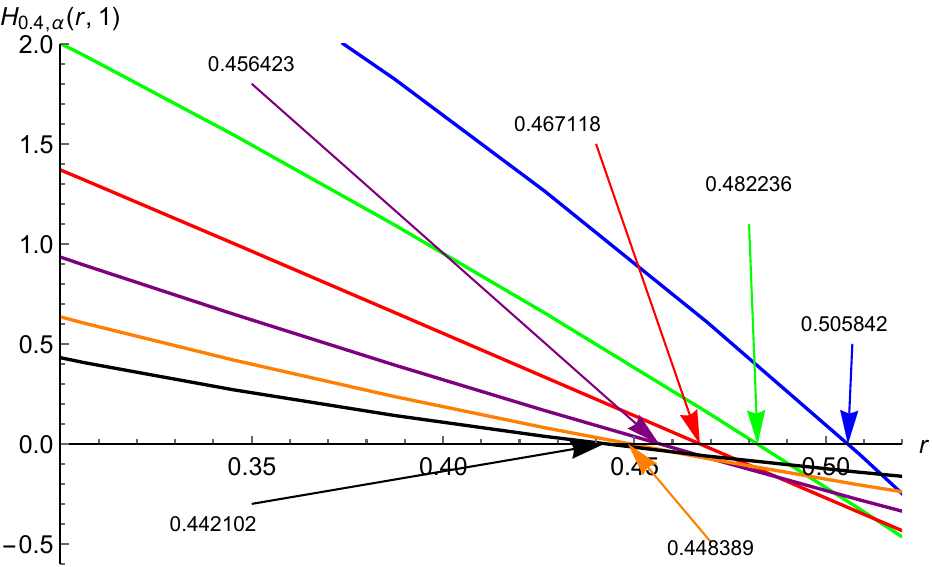}
	\end{center}
	\caption{The graph of $H_{0.1, \alpha}(r,1)$ and $H_{0.4, \alpha}(r,1)$ when $\alpha=0.5,1.0,1.5,2.0,2.5, 3.0$.}
	\label{figure-4.3-a}
\end{figure}

\begin{table}[ht]
	\centering
	\begin{tabular}{|l|l|l|}
		\hline
		$\alpha$& $r_{0.7}(\alpha,1)$& $r_{0.9}(\alpha,1)$ \\
		\hline
		$(0,0.5]$& $(0.551329 \searrow 0.520253]$& $(0.551329 \searrow 0.528847]$\\
		\hline
		$(0.5,1]$& $(0.520253 \searrow 0.505149]$& $(0.528847 \searrow 0.518328]$\\
		\hline
		$(1,1.5]$& $(0.505149 \searrow 0.495956]$& $(0.518328 \searrow 0.512104]$\\
		\hline
		$(1.5,2]$& $(0.495956 \searrow 0.489709]$& $(0.512104 \searrow 0.507962]$\\
		\hline
		$(2,2.5]$& $(0.489709 \searrow 0.485167]$& $(0.507962 \searrow 0.504999]$\\
		\hline
		$(2.5,3]$& $(0.485167 \searrow 0.481706]$& $(0.504999 \searrow 0.502771]$\\
		\hline
	\end{tabular}
	\vspace{3mm}
	\caption{Values of $r_{0.7}(\alpha,1)$ and $r_{0.9}(\alpha,1)$ for various values of $\alpha$.}
	\label{tabel-4.3-b}
\end{table}

\begin{figure}[!htb]
	\begin{center}
		\includegraphics[width=0.48\linewidth]{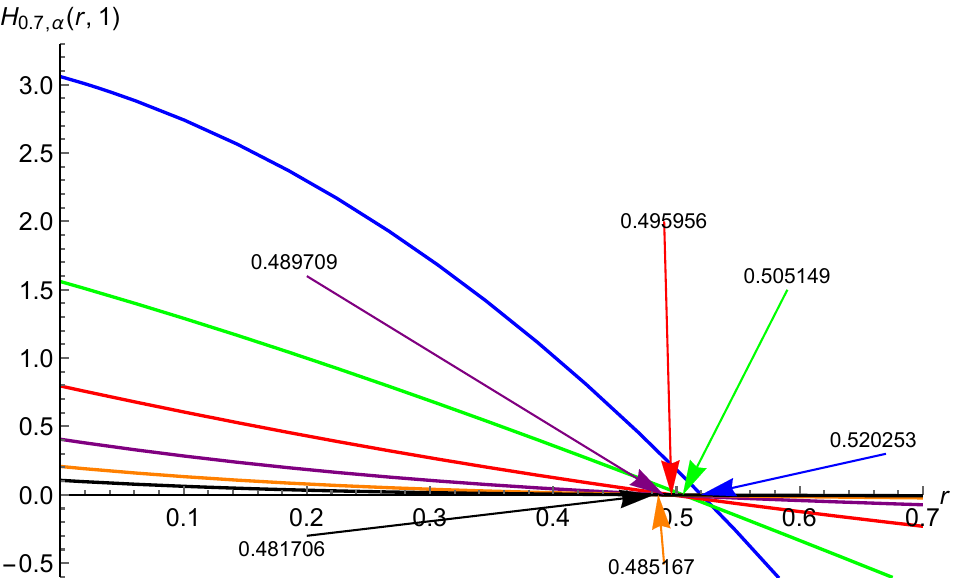}
		\,
		\includegraphics[width=0.48\linewidth]{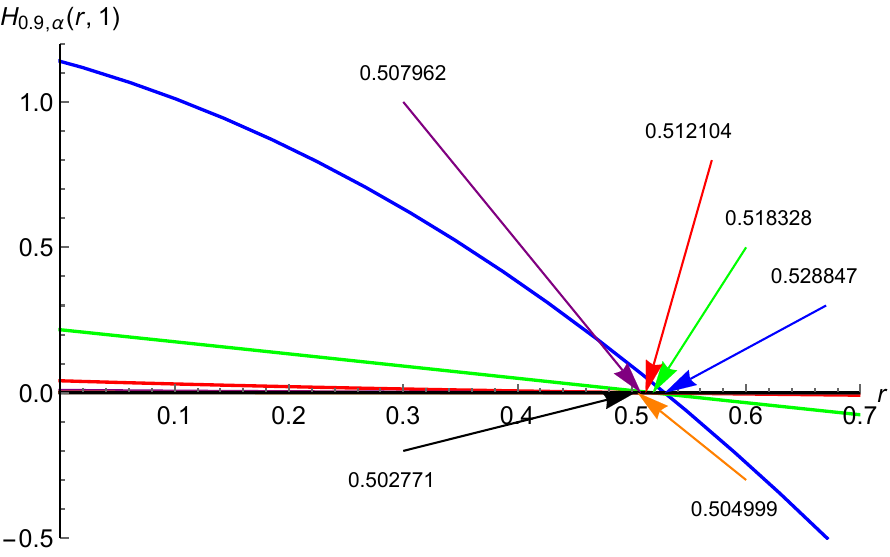}
	\end{center}
	\caption{The graph of $H_{0.7, \alpha}(r,1)$ and $H_{0.9, \alpha}(r,1)$ when $\alpha=0.5,1.0,1.5,2.0,2.5, 3.0$.}
	\label{figure-4.3-b}
\end{figure}
By adopting the similar arguments as in Corollary \ref{him-vasu-P5-cor-4.15} and using Comparison Principle \ref{him-vasu-P5-thm-2.1-comparison-principle}, we obtain the following result.
\begin{cor}
	Let $\Omega_{1}, \Omega_{2}, \ldots , \Omega_{n}$ be finite $n$- proper simply connected domains in $\mathbb{C}$ containing $\mathbb{D}$ such that $\Omega_{1} \subseteq \Omega_{2} \subseteq \ldots \subseteq \Omega_{n}$. Suppose that Bloch-Bohr radius exists for all $\mathcal{B}_{\mathcal{H},\Omega_{1}}(\alpha), \mathcal{B}_{\mathcal{H},\Omega_{2}}(\alpha), \ldots, $  $\mathcal{B}_{\mathcal{H},\Omega_{n}}(\alpha)$, and call them as $r_{1}, r_{2}, \ldots,  r_{n}$ respectively. Then $r_{1} \leq r_{2} \ldots \leq r_{n}$.
\end{cor}

Now we see that the $p$-Bohr radius for functions in $\mathcal{B}^{*}_{\mathcal{H},\Omega_{\gamma}}(\alpha)$ is $0$. To see this, we consider the function 
\begin{equation*}
f_{\gamma}(z)= h(z) + \overline{g(z)}= \frac{1}{(1-\gamma)(1-z)} + \overline{\frac{z}{(1-\gamma)(1-z)}}, \,\,\,\, z \in \Omega_{\gamma}
\end{equation*}
so that $f_{\gamma}(z)=\sum_{n=0}^{\infty} a_{n} z^n + \overline{\sum_{n=1}^{\infty} b_{n} z^n}$ in $\mathbb{D}$, where
$a_{0}=1/(1-\gamma)$, $a_{n}=b_{n}=1/(1-\gamma)$ for $n \geq 1$. It is easy to see that $f_{\gamma} \in \mathcal{B}^{*}_{\mathcal{H},\Omega_{\gamma}}(\alpha)$ and $||f||^{*}_{\mathcal{H},\Omega_{\gamma},\alpha}=|a_{0}|=1/(1-\gamma)$ for $\alpha>0$. Therefore, we obtain
\begin{equation*}
|a_{0}| + \sum \limits_{n=1}^{\infty}\left(|a_{n}|^p + |b_{n}|^p\right)^{1/p}\, r^n = \frac{1}{1-\gamma} + \sum \limits_{n=1}^{\infty} \frac{2^{1/p}}{1-\gamma}\, r^n>\frac{1}{1-\gamma}=||f||^{*}_{\mathcal{H},\Omega_{\gamma},\alpha}
\end{equation*}
for all $r \in (0,1)$, which shows that the $p$-Bohr radius for $f$ is $0$. But, if we consider an additional condition, namely sense-preserving to the functions in $\mathcal{B}^{*}_{\mathcal{H},\Omega_{\gamma}}(\alpha)$, then the $p$-Bohr radius exists and is obtained in the following result.

\begin{thm} \label{him-vasu-P5-thm-4.4}
Let $f=h+\overline{g} \in \mathcal{B}^{*}_{\mathcal{H}, \Omega_{\gamma}}(\alpha)$ be a sense-preserving harmonic mapping such that $||f||^{*}_{\mathcal{H}, \Omega_{\gamma}} (\alpha)\leq 1$. If $|g'(z)| \leq k |h'(z)|$ in $\mathbb{D}$ and 
\begin{equation*}
h(z)=\sum \limits_{n=1}^{\infty}  a_{n} z^{n} \,\,\, \mbox{and} \,\,\, g(z)=\sum \limits_{n=1}^{\infty}  b_{n} z^{n} \,\,\, \mbox{in} \,\, \mathbb{D},
\end{equation*}
then, for each $p \geq 1$, we have 
\begin{equation*}
|a_{0}| + \sum \limits_{n=1}^{\infty} \left(|a_{n}|^p + |b_{n}|^p\right)^{1/p}\, r^n \leq 1
\end{equation*}
for $|z|=r \leq r^{*}_{\gamma}(\alpha,p,d)$, where $r^{*}_{\gamma}(\alpha,p,d)$ is the smallest root of $H_{2}(r)=0$ in $(0,1)$, where
\begin{equation*}
H_{2}(r)=6(1-((1-\gamma)r + \gamma)^2)^{2\alpha} \left((1+dr)^2 - k^2 \, (r+d)^2\right)- K_{p}\,(1+k^2)\pi ^2 \,(1+d)^2 \, (1-\gamma)^{2 \alpha} r^2,
\end{equation*}
$d = |g'(0)|/(k\, |h'(0)|)$, and $K_{p}=\max \{2^{(2/p)-1}, 1\}$.
\end{thm}

\begin{pf}
For the sense-preserving harmonic mapping $f=h+\overline{g}$ in $\Omega_{\gamma}$ satisfying $|g'(z)| \leq k |h'(z)|$ in $\mathbb{D}$, the dilation $\omega=\omega_{f}=g'/h'$ satisfies $|\omega_{f}(z)| \leq k<1$ for $z \in \mathbb{D}$. If $\omega_{f}$ is non-constant, then by the maximum modulus principle, $|\omega_{f}(z)| <k$ for all $z \in \mathbb{D}$. Thus, if we assume that $\omega_{f}$ is non-constant then there exists an analytic function $\psi _{\omega} : \mathbb{D} \rightarrow \mathbb{D}$ such that $\psi_{\omega}=\omega_{f}/k$ in $\mathbb{D}$ {\it i.e.} $g'(z)=k\psi_{\omega}(z)h'(z)$ for $z \in \mathbb{D}$. Therefore, 
\begin{equation*}
 J_{f}(z)=|h'(z)|^2 (1-|\omega(z)|^2)\,\, {\it i.e.,}\,\, |h'(z)|=\sqrt{\frac{|J_{f}(z)|}{1-|\omega(z)|^2}}= \sqrt{\frac{|J_{f}(z)|}{1-k^2\,|\psi_{\omega}(z)|^2}}, \,\,\, z \in \mathbb{D}.
\end{equation*}
Since $|\psi_{\omega}(z)|<1$ in $\mathbb{D}$, from the Pick's conformally invariant form of the Schwarz lemma, we have 
\begin{equation} \label{him-vasu-P5-e-4.19}
|\psi_{\omega}(z)| \leq \frac{|z| + |\psi_{\omega}(0)|}{1+|\psi_{\omega}(0)|\, |z|}, \,\,\, z \in \mathbb{D}.
\end{equation}	
The given assumption $||f||^{*}_{\mathcal{H}, \Omega_{\gamma}} (\alpha)\leq 1$ gives 
\begin{equation*}
\frac{\sqrt{J_{f}(z)}}{\lambda ^{\alpha}_{\Omega_{\gamma}}(z)} \leq 1- |a_{0}|,
\end{equation*}	
which is equivalent to 
\begin{equation} \label{him-vasu-P5-e-4.20}
\left(1- ((1-\gamma)|z|+\gamma)^2\right)^{\alpha} \, \sqrt{J_{f}(z)} \leq (1-\gamma)^{\alpha}\, (1-|a_{0}|), \,\,\,\, z \in \Omega_{\gamma}.
\end{equation}
Since \eqref{him-vasu-P5-e-4.20} holds for all $z \in \Omega_{\gamma}$ and $\mathbb{D} \subseteq \Omega_{\gamma}$, then \eqref{him-vasu-P5-e-4.20} also holds for all $z \in \mathbb{D}$.
Thus, from \eqref{him-vasu-P5-e-4.19} and \eqref{him-vasu-P5-e-4.20}, we obtain
\begin{align} \label{him-vasu-P5-e-4.21}
|h'(z)|&=\sqrt{\frac{|J_{f}(z)|}{1-k^2 \,|\psi_{\omega}(z)|^2}}\\ \nonumber
&\leq \frac{(1-\gamma)^{\alpha}\, (1-|a_{0}|)}{\left(1- ((1-\gamma)|z|+\gamma)^2\right)^{\alpha}}\,\,\frac{1}{\sqrt{1-k^2 \,|\psi_{\omega}(z)|^2}}\\ \nonumber
& \leq \frac{(1-\gamma)^{\alpha}\, (1-|a_{0}|)}{\left(1- ((1-\gamma)|z|+\gamma)^2\right)^{\alpha}} \,\, \left(1-k^2 \, \left(\frac{|z| + |\psi_{\omega}(0)|}{1+|\psi_{\omega}(0)|\, |z|}\right)^2\right)^{\frac{-1}{2}}\\ \nonumber
&=  \frac{(1-\gamma)^{\alpha}\, (1-|a_{0}|)}{\left(1- ((1-\gamma)|z|+\gamma)^2\right)^{\alpha}}  \,\, \frac{1+|\psi_{\omega}(0)||z|}{\sqrt{(1+|\psi_{\omega}(0)| |z|)^2 - k^2 (|z| + |\psi_{\omega}(0)|)^2}}\\ \nonumber
&\leq  \frac{(1-\gamma)^{\alpha}\, (1-|a_{0}|)}{\left(1- ((1-\gamma)|z|+\gamma)^2\right)^{\alpha}}  \,\, \frac{1+|\psi_{\omega}(0)|}{\sqrt{(1+|\psi_{\omega}(0)| |z|)^2 - k^2 (|z| + |\psi_{\omega}(0)|)^2}}, \,\,\, z \in \mathbb{D}.
\end{align}
Since $|g'(z)| \leq k\, |h'(z)|$ in $\mathbb{D}$ then $ |h'(z)|^2 + |g'(z)|^2 \leq  (1+k^2)|h'(z)|^2$. Using the fact $|h'(z)|^2 + |g'(z)|^2 \leq (|h'(z)| + |g'(z)|)^2$, the inequality \eqref{him-vasu-P5-e-4.21} leads to
\begin{equation} \label{him-vasu-P5-e-4.22}
|h'(z)|^2 + |g'(z)|^2 \leq \frac{(1+k^2)(1-\gamma)^{2\alpha}\, (1-|a_{0}|)^2}{\left(1- ((1-\gamma)|z|+\gamma)^2\right)^{2\alpha}}  \,\, \frac{(1+|\psi_{\omega}(0)|)^2}{(1+|\psi_{\omega}(0)| |z|)^2 - k^2 (|z| + |\psi_{\omega}(0)|)^2}
\end{equation}
in $\mathbb{D}$. Integrating \eqref{him-vasu-P5-e-4.22} over the circle $|z|=r<1$, we obtain 
\begin{align} \label{him-vasu-P5-e-4.23}
&\sum \limits_{n=1}^{\infty} n^2 (|a_{n}|^2 + |b_{n}|^2)r^{2(n-1)} \\ \nonumber
&\leq \frac{(1+k^2)(1-\gamma)^{2\alpha}\, (1-|a_{0}|)^2}{\left(1- ((1-\gamma)r+\gamma)^2\right)^{2\alpha}}  \,\, \frac{(1+|\psi_{\omega}(0)|)^2}{(1+|\psi_{\omega}(0)| r)^2 - k^2 (r + |\psi_{\omega}(0)|)^2}.
\end{align}
Using the Cauchy-Schwarz inequality and  \eqref{him-vasu-P5-e-4.23}, we obtain
\begin{align} \label{him-vasu-P5-e-4.24}
&|a_{0}| + \sum \limits_{n=1}^{\infty}\left(|a_{n}|^p + |b_{n}|^p\right)^{1/p}\, r^n \\ \nonumber
 & \leq 
|a_{0}| + \sqrt{\sum \limits_{n=1}^{\infty} n^2 (|a_{n}|^p + |b_{n}|^p) ^{2/p} \, r^{2n}} \,\, \sqrt{\sum \limits_{n=1}^{\infty} \frac{1}{n^2}} \\ \nonumber
& \leq |a_{0}| + \sqrt{K_{p} \, \sum \limits_{n=1}^{\infty} n^2 (|a_{n}|^2 + |b_{n}|^2) \, r^{2n}}\,\,\, \sqrt{\frac{\pi ^2}{6}} \\ \nonumber
& \leq |a_{0}| + \frac{\sqrt{K_{p}(1+k^2)}\,(1-\gamma)^{\alpha}\, (1-|a_{0}|)}{\left(1- ((1-\gamma)r+\gamma)^2\right)^{\alpha}}  \,\, \frac{(1+|\psi_{\omega}(0)|)r}{\sqrt{(1+|\psi_{\omega}(0)| r)^2 - k^2 (r + |\psi_{\omega}(0)|)^2}}\,\, \sqrt{\frac{\pi ^2}{6}} \\ \nonumber
&\leq 1,
\end{align} 
provided $r \leq r^{*}_{\gamma}(\alpha,p)$, where $r^{*}_{\gamma}(\alpha,p)$ is the smallest root of 
\begin{equation*}
\frac{\sqrt{K_{p}(1+k^2)}\,(1-\gamma)^{\alpha}}{\left[1- ((1-\gamma)r+\gamma)^2\right]^{\alpha}}  \,\, \frac{(1+|\psi_{\omega}(0)|)r}{\sqrt{(1+|\psi_{\omega}(0)| r)^2 - k^2 (r + |\psi_{\omega}(0)|)^2}}\,\, \sqrt{\frac{\pi ^2}{6}}=1, 
\end{equation*}
or equivalently,
\begin{align} \label{him-vasu-P5-e-4.25}
&
6(1-((1-\gamma)r + \gamma)^2)^{2\alpha} \left((1+dr)^2  - k^2 \, (r+d)^2\right) \\ \nonumber 
&-K_{p}\,(1+k^2)\pi ^2 \,(1+d)^2 \, (1-\gamma)^{2 \alpha} r^2=0
\end{align}
in $(0,1)$, where $d = |\psi_{\omega}(0)|=|g'(0)|/(k\, |h'(0)|)$. Using the similar lines of argument as in the proof of Theorem \ref{him-vasu-P5-thm-4.2}, we can show that \eqref{him-vasu-P5-e-4.25} has root in $(0,1)$ and choose $r^{*}_{\gamma}(\alpha,p)$ to be the smallest root in $(0,1)$. This completes the proof.
\end{pf}

It is worth to point out that $k \rightarrow 1$ corresponds to the sense-preserving harmonic mappings. Hence, by taking $k \rightarrow 1$ in Theorem \ref{him-vasu-P5-thm-4.4}, we obtain the following corollary.

\begin{cor} \label{him-vasu-P5-cor-4.26}
Let $f=h+\overline{g} \in \mathcal{B}^{*}_{\mathcal{H}, \Omega_{\gamma}}(\alpha)$ with $||f||^{*}_{\mathcal{H}, \Omega_{\gamma}} (\alpha)\leq 1$. Let $h(z)=\sum \limits_{n=1}^{\infty}  a_{n} z^{n}$, $g(z)=\sum \limits_{n=1}^{\infty}  b_{n} z^{n}$ in $\mathbb{D}$ and $f=h+\overline{g}$ is sense-preserving in $\mathbb{D}$.
Then, for each $p \geq 1$, we have 
\begin{equation*}
|a_{0}| + \sum \limits_{n=1}^{\infty} \left(|a_{n}|^p + |b_{n}|^p\right)^{1/p}\, r^n \leq 1
\end{equation*}
for $|z|=r \leq r^{*}_{\gamma}(\alpha,p,d_{1})$, where $r^{*}_{\gamma}(\alpha,p,d_{1})$ is the smallest root of $H_{3}(r)=0$ in $(0,1)$, where
\begin{equation*}
H_{3}(r):=6(1-((1-\gamma)r + \gamma)^2)^{2\alpha} \left((1+d_{1}r)^2 -  \, (r+d_{1})^2\right)- 2\,K_{p}\,\pi ^2 \,(1+d_{1})^2 \, (1-\gamma)^{2 \alpha} r^2
\end{equation*}
and $d_{1} = |g'(0)|/|h'(0)|$ and $K_{p}=\max \{2^{(2/p)-1}, 1\}$.
\end{cor}

\begin{figure}[!htb]
	\begin{center}
		\includegraphics[width=0.48\linewidth]{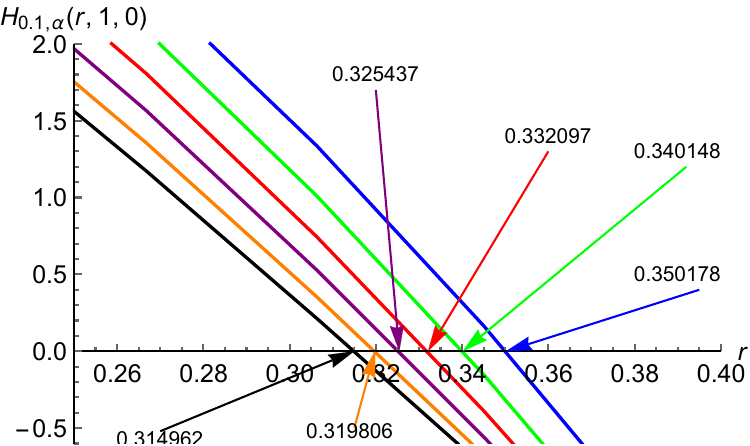}
		\,
		\includegraphics[width=0.48\linewidth]{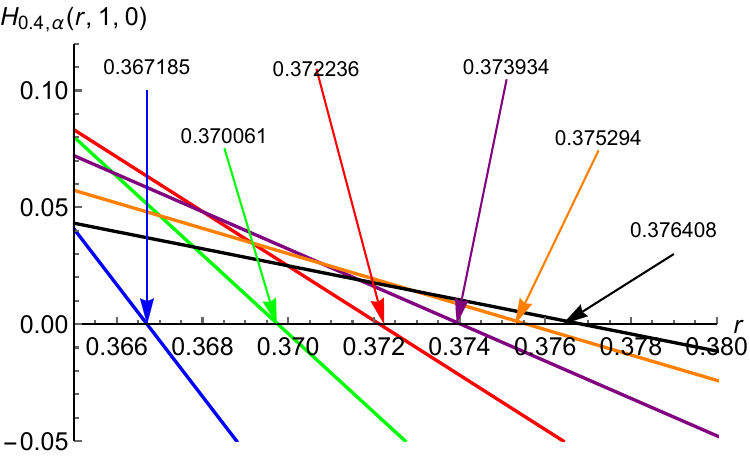}
	\end{center}
	\caption{The graph of $H_{0.1, \alpha}(r,1,0)$ and $H_{0.4, \alpha}(r,1,0)$ when $\alpha=0.5,1.0,1.5,2.0,2.5, 3.0$.}
	\label{figure-4.4-a}
\end{figure}

\begin{figure}[!htb]
	\begin{center}
		\includegraphics[width=0.48\linewidth]{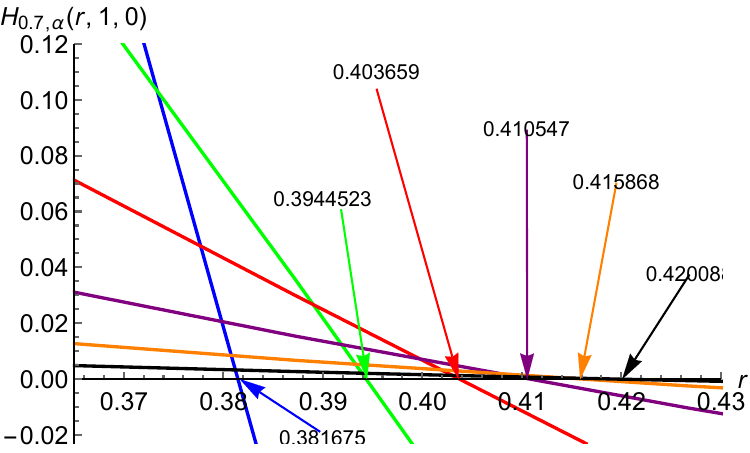}
		\,
		\includegraphics[width=0.48\linewidth]{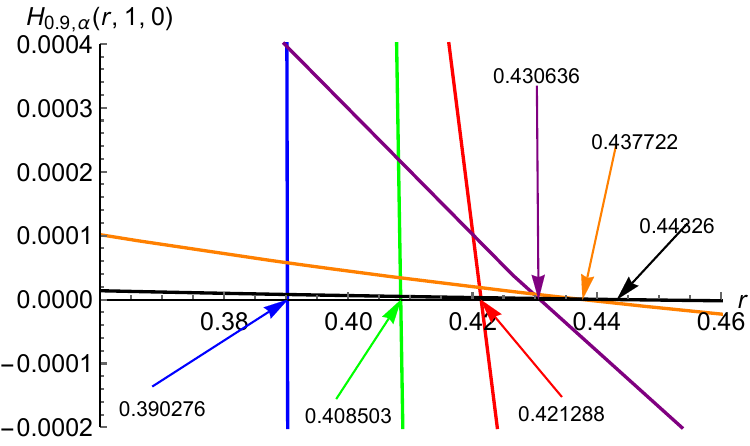}
	\end{center}
	\caption{The graph of $H_{0.7, \alpha}(r,1,0)$ and $H_{0.9, \alpha}(r,1,0)$ when $\alpha=0.5,1.0,1.5,2.0,2.5, 3.0$.}
	\label{figure-4.4-b}
\end{figure}

\begin{table}[ht]
	\centering
	\begin{tabular}{|l|l|l|l|l|}
		\hline
		$\alpha$& $r^{*}_{0.1}(\alpha,1,0)$& $r^{*}_{0.4}(\alpha,1,0)$ & $r^{*}_{0.7}(\alpha,1,0)$ &$r^{*}_{0.9}(\alpha,1,0)$ \\
		\hline
		$0.0$& $0.363223$& $0.363223$ & $0.363223$ &$0.363223$\\
		\hline
		$0.5$& $0.350178$& $0.367185$ &$0.381675$ &$0.390276$\\
		\hline
		$1.0$& $0.340148$& $0.370061$ & $0.394452$ &$0.408503$\\
		\hline
		$2.0$& $0.325437$& $0.373934$ &$0.410547$ &$0.430636$\\
		\hline
		$2.5$& $0.319806$& $0.375294$ &$0.415868$ &$0.437722$\\
		\hline
		$3.0$& $0.314962$& $0.376408$ &$0.420088$ &$0.443260$\\
		\hline
	\end{tabular}
	\vspace{3mm}
	\caption{Values of $r^{*}_{\gamma}(\alpha, p, d_{1})$ for different values of  $\gamma, \alpha$ when $p=1$, $d_{1}=0$.} 
	\label{tabel-4.3-cor}
\end{table}
\noindent Denote $H_{3}= H_{\gamma,\alpha}(r,p,d_{1})$. In Table \ref{tabel-4.3-cor}, for $p=1$ and $d_{1}=0$, we compute $p$- Bloch-Bohr radius $r^{*}_{\gamma}(\alpha,p,d_{1})$ for $f \in \mathcal{B}^{*}_{\mathcal{H}, \Omega_{\gamma}}(\alpha) $ which is sense-preserving in $\mathbb{D}$. From Table \ref{tabel-4.3-cor}, we observe that $r^{*}_{\alpha, p, d_{1}}$ is increasing in $\gamma$ for fixed values of $\alpha$ and $p$. \\[2mm]

By using the similar arguments as in Corollary \ref{him-vasu-P5-cor-4.15} and using the Comparison Principle \ref{him-vasu-P5-thm-2.1-comparison-principle}, we obtain the following Corollary.
\begin{cor}
	Let $\Omega_{1}, \Omega_{2}, \ldots , \Omega_{n}$ be finite $n$-proper simply connected domains in $\mathbb{C}$ containing $\mathbb{D}$ such that $\Omega_{1} \subseteq \Omega_{2} \subseteq \ldots \subseteq \Omega_{n}$. Suppose that Bloch-Bohr radius exists for all $\mathcal{B}^{*}_{\mathcal{H},\Omega_{1}}(\alpha), \mathcal{B}^{*}_{\mathcal{H},\Omega_{2}}(\alpha), \ldots, $  $\mathcal{B}^{*}_{\mathcal{H},\Omega_{n}}(\alpha)$, and call them as $r^{*}_{1}, r^{*}_{2},\ldots , r^{*}_{n}$ respectively. Then $r^{*}_{1} \leq r^{*}_{2} \leq \ldots \leq r^{*}_{n}$.
\end{cor}

We note that $K_{p}=\max \{2^{(2/p)-1}, 1\}=1$ when $p>1$. Then, for $p>1$, we have the following result.

\begin{cor} \label{him-vasu-P5-cor-4.26-a}
	Let $f=h+\overline{g} \in \mathcal{B}^{*}_{\mathcal{H}, \Omega_{\gamma}}(\alpha)$ with $||f||^{*}_{\mathcal{H}, \Omega_{\gamma}} (\alpha)\leq 1$. If $h(z)=\sum \limits_{n=1}^{\infty}  a_{n} z^{n}$, $g(z)=\sum \limits_{n=1}^{\infty}  b_{n} z^{n}$ in $\mathbb{D}$ and $f=h+\overline{g}$ is sense-preserving in $\mathbb{D}$
	then, for each $p \geq 2$, we have 
	\begin{equation*}
	|a_{0}| + \sum \limits_{n=1}^{\infty} \left(|a_{n}|^p + |b_{n}|^p\right)^{1/p}\, r^n \leq 1
	\end{equation*}
	for $|z|=r \leq r^{*}_{\gamma}(\alpha,p,d_{1})$, where $r^{*}_{\gamma}(\alpha,p,d_{1})$ is the smallest root of $H_{4}(r)=0$ in $(0,1)$, where
	\begin{equation*}
	H_{4}(r):=6(1-((1-\gamma)r + \gamma)^2)^{2\alpha} \left((1+d_{1}r)^2 -  \, (r+d_{1})^2\right)- 2\,\pi ^2 \,(1+d_{1})^2 \, (1-\gamma)^{2 \alpha} r^2
	\end{equation*}
	and $d_{1} = |g'(0)|/|h'(0)|$ and $K_{p}=\max \{2^{(2/p)-1}, 1\}$.
\end{cor}

An observation of Corollary \ref{him-vasu-P5-cor-4.26-a} shows that the $p$-Bohr radius $r^{*}_{\gamma}(\alpha,p,d_{1})$ depends on $d_{1}$ {\it i.e.,} on $|\omega_{f}(0)|$, which can be seen from the following example. Liu and Ponnusamy \cite{Liu-Results-Math-2018} have considered the following example and have shown the dependence of $|\omega_{f}(0)|$ about $p$-Bohr radius. 

\begin{example}
Let $\gamma=0$ {\it i.e.,} $\Omega_{\gamma}=\mathbb{D}$. Fix $\lambda \in [1/2, 1)$. Let $F_{\lambda}$ be the following one parameter family of functions (see \cite[Example 3]{Liu-Results-Math-2018}) defined by 
$$
F_{\lambda}(z)= H_{\lambda} (z) + \overline{ G_{\lambda} (z)} \,\,\,\,\,\,\,\,\, \mbox{for} \,\,\, z \in \mathbb{D},
$$
where 
$$
H_{\lambda} (z)= 1-2\sqrt{\lambda-\lambda^2} + \frac{1}{2}\, \log \frac{1+z}{1-z}
$$
and
$$
 G_{\lambda} (z)= \frac{\lambda-1}{2}\, \log(1-z^2) + \frac{\lambda}{2}\, \log \frac{1+z}{1-z}.
$$
Furthermore, let 
$$
F_{\lambda}(z)= \sum \limits_{n=0}^{\infty} a_{n,\lambda}\, z^n + \overline{\sum \limits_{n=1}^{\infty} b_{n,\lambda}\, z^n} \,\,\,\,\,\,\,\,\, \mbox{for} \,\, z \in \mathbb{D}.
$$
A simple computation shows that 
\begin{equation*}
H'_{\lambda} (z)= \frac{1}{1-z^2} \,\,\, \mbox{and} \,\,\, 
G'_{\lambda}(z)=\frac{(1-\lambda)z+\lambda}{1-z^2} \,\, \,\,\, \mbox{for} \,\,\,  z \in \mathbb{D}.
\end{equation*}
Then the dilation $\omega_{F_{\lambda}}(z)=(1-\lambda)z +\lambda$ and $|\omega_{F_{\lambda}}(z)|<1$ in $\mathbb{D}$. 
Hence, $F_{\lambda}$ is sense-preserving in $\mathbb{D}$. Now, we see that 
\begin{equation} \label{him-vasu-P5-e-4.27}
|G'_{\lambda}(z)| \geq \frac{\lambda-(1-\lambda)|z|}{|1-z^2|}\geq \frac{2\lambda -1}{|1-z^2|} \,\, \,\,\, \mbox{for} \,\,\,  z \in \mathbb{D}.
\end{equation} 
Thus, from \eqref{him-vasu-P5-e-4.27}, we obtain 
\begin{equation*}
(1-|z|^2)\, \sqrt{J_{F_{\lambda}}(z)} \leq (1-|z|^2)\, \sqrt{\frac{1}{|1-z^2|^2} - \frac{(2\lambda -1)^2}{|1-z^2|^2}} \leq 2\sqrt{\lambda - \lambda ^2}, \,\,\, z \in \mathbb{D},
\end{equation*}
which shows that $F_{\lambda} \in \mathcal{B}^{*}_{\mathcal{H}}(1)$. Furthermore, for $x \in (-1,0)$, we have 
\begin{equation*}
(1-|x|^2)\, \sqrt{J_{F_{\lambda}}(x)} \rightarrow 2 \sqrt{\lambda - \lambda ^2}
\end{equation*}
as $x \rightarrow -1 ^{+}$, which implies that $\beta ^{*}_{\mathcal{H}}(\alpha)= 2 \sqrt{\lambda - \lambda ^2}$ and hence, $||F_{\lambda}||_{\mathcal{B}^{*}_{\mathcal{H}}(1)}=1$.
A simple computation shows that 
\begin{equation*}
|a_{0,\lambda}|+ \sum \limits_{n=1}^{\infty} \left(|a_{n,\lambda}|^p + |b_{n,\lambda}|^p\right)^{1/p} \,\, r^{n} > |a_{0,\lambda}|=1- 2 \sqrt{\lambda - \lambda ^2}
\end{equation*}
for each $r \in (0,1)$. Now, we observe that $|\omega_{F_{\lambda}}(0)|= |G'_{\lambda}(0)| /|H'_{\lambda}(0)|=\lambda$. By taking $\lambda \rightarrow 1^{-}$, we see that $1- 2 \sqrt{\lambda - \lambda ^2} \rightarrow 1=||F_{\lambda}||_{\mathcal{B}^{*}_{\mathcal{H}}(1)}$, which shows that $p$-Bohr radius for $F_{\lambda}$ tends to $0$ when $\lambda \rightarrow 1^{-}$. Therefore, the $p$-Bohr radius depends on $|\omega_{F_{\lambda}}(0)|$.
\end{example} 

\noindent\textbf{Acknowledgment:} The first author is supported by SERB-CORE Grant and second author is supported by CSIR (File No: 09/1059(0020)/2018-EMR-I), New Delhi, India.\\

\noindent\textbf{Competing interests declaration:} The authors declare none.


\begin{thebibliography}{99}
	
	\bibitem{abdulhadi-20} {\sc Z. Abdulhadi} and  {\sc Y. Abu-Muhanna}, Landau’s theorem for biharmonic mappings, {\it J. Math. Anal. Appl.} {\bf  338}  (2008), 705--709.
	
	\bibitem{Abu-2010} {\sc Y. Abu-Muhanna},  Bohr's phenomenon in subordination and bounded harmonic classes, {\it Complex Var. Elliptic Equ.} {\bf  55} (2010), 1071--1078.
	
	
	\bibitem{abu-2011} {\sc Y. Abu-Muhanna} and  {\sc R. M. Ali}, Bohr's phenomenon for analytic functions into the exterior of a compact convex body, {\it J. Math. Anal. Appl.} {\bf  379}  (2011), 512--517.
	
	
	
	
	\bibitem{abu-2014} {\sc Y. Abu Muhanna, R. M. Ali, Z. C. Ng,} and  {\sc S. F. M Hasni}, Bohr radius for subordinating families of analytic functions and bounded harmonic mappings, 
	{\it J. Math. Anal. Appl.} {\bf 420} (2014), 124--136.
	
	
	\bibitem{Ahamed-Allu-Halder-P3-2020} {\sc M. B. Ahamed, V. Allu} and {\sc H. Halder}, The Bohr Phenomenon for analytic functions on simply connected domains, {\it  Ann. Acad. Sci. Fenn. Ser. A I Math.} (2021), To appear.
	
	\bibitem{Ahamed-Allu-Halder-AAMP-2020} {\sc M. B. Ahamed, V. Allu} and {\sc H. Halder}, Bohr radius for certain classes of close-to-convex harmonic mappings, {\it Anal. Math. Phys.} {\bf 11} (111) (2021).
	
	
	\bibitem{aizn-2000} {\sc L. Aizenberg}, Multidimensional analogues of Bohr's theorem on power series, \textit{Proc. Amer. Math. Soc.} {\bf 128} (2000), 1147--1155.
	
	
	\bibitem{aizenberg-2001} {\sc L. Aizenberg, A. Aytuna}  and {\sc P. Djakov}, Generalization of theorem on Bohr for bases in spaces of holomorphic functions of several complex variables, 
	{\it J. Math. Anal. Appl.} {\bf  258} (2001), 429--447.
	
	
	\bibitem{aizn-2007} {\sc L. Aizenberg}, Generalization of results about the Bohr radius for power series, {\it Stud. Math.}  {\bf 180}  (2007), 161--168.  
	
	\bibitem{aizenberg-2012} {\sc L. Aizenebrg}, Remarks on the Bohr and Rogosinski phenomena for power series, {\it Anal. Math. Phys.} {\bf 2} (2012), 69--78.
	
	
	
	
	
	
	\bibitem{Ali-2017} {\sc R. M. Ali, R.W. Barnard} and {\sc  A.Yu. Solynin}, A note on Bohr's phenomenon for power series, {\it J. Math. Anal. Appl.} {\bf 449} (2017), 154-167.
	
	
	\bibitem{alkhaleefah-2019} {\sc S. A. Alkhaleefah, I.R. Kayumov} and {\sc S. Ponnusamy}, On the Bohr inequality with a fixed zero coefficient, {\it Proc. Amer. Math. Soc.} {\bf 147} (2019), 5263--5274.
	
	\bibitem{Himadri-Vasu-P1} {\sc V. Allu} and {\sc H. Halder},  Bhor phenomenon for certain subclasses of Harmonic Mappings, {\it Bull. Sci. Math.} {\bf 173} (2021), 103053.
	
	\bibitem{Himadri-Vasu-P2} {\sc V. Allu} and {\sc H. Halder}, Bohr radius for certain classes of starlike and convex univalent functions, {\it J. Math. Anal. Appl.} {\bf 493}(1) (2021), 124519.
	
		\bibitem{Himadri-Vasu-P3} {\sc V. Allu} and {\sc H. Halder}, Bohr phenomenon for certain close-to-convex analytic functions, {\it Comput. Methods Funct. Theory} (2021), To appear.
	
	
	\bibitem{Anderson-1974} {\sc J. M. Anderson, J. Clunie,} and {\sc Ch. Pommerenke}, On Bloch functions and normal functions, {\it J. Reine. Anjew. Math.} {\bf 270} (1974), 12--37.
	
	\bibitem{Arazy-1985} {\sc J. Arazy, S. D. Fisher,} and {\sc J. Peetre}, M\"{o}bius invariant function spaces, {\it J. Reine. Anjew. Math.} {\bf 363} (1985), 110--145.
	
	
	
	\bibitem{Ayt & Dja & BLMS & 2013} {\sc A. Aytuna} and {\sc P. Djakov}, Bohr property of bases in the space of entire functions and its generalizations, {\it Bull. London Math. Soc.} \textbf{45}(2)(2013), 411--420.
	
	\bibitem{Bai-CAOT-2019} {\sc X.-X. Bai, D. Pellegrino}, and {\sc J. B. Seoane-Sep$\rm \acute{U}$lveda},  Landau-type theorems of poly-harmonic mappings and log-$p$-harmonic mappings, {\it Complex Anal. Oper. Theory} {\bf 13} (2019), 321--340.
	
	\bibitem{bayart-advance-2014} {\sc F. Bayart, D. Pellegrino}, and {\sc J. B. Seoane-Sep$\rm \acute{U}$lveda}, The Bohr radius of the $n$-dimensional polydisk is equivalent to $\sqrt{(log \, n)/n}$, {\it Adv. Math.} {\bf 264} (2014), 726--746.
	
	\bibitem{beardon-minda-hyperbolic-density} {\sc A. F. Beardon} and {\sc D. Minda}, The hyperbolic metric and geometric function theory, {\it Proceedings of the International Workshop on Quasiconformal Mappings and their Applications (IWQCMA05)}.
	
	\bibitem{bene-2004} {\sc C. B$ {\rm \acute{E}} $n$ {\rm \acute{E}} $teau}, {\sc A. Dahlner} and {\sc D. Khavinson}, Remarks on the Bohr phenomenon, {\it Comput. Methods Funct. Theory} \textbf{4}(1) (2004), 1-19.
	
	\bibitem{bhowmik-2021} {\sc B. Bhowmik} and {\sc N. Das}, Bohr phenomenon for operator-valued functions, {\it 
		Proc. Edinburgh Math. Soc.}, https://doi.org/10.1017/S0013091520000395.
	
	
	\bibitem{boas-1997} {\sc H.P. Boas} and {\sc D. Khavinson}, Bohr's power series theorem in several variables, {\it Proc. Amer. Math. Soc}  {\bf 125} (1997), 2975--2979.
	
	
	\bibitem{Bohr-1914} {\sc H. Bohr}, A theorem concerning power series,  {\it Proc. Lond. Math. Soc}. s2-13 (1914), 1--5.
	
	\bibitem{Bonk-CMFT-1994} {\sc M. Bonk, D. Minda,} and {\sc H. Yanagihara}, The hyperbolic metric on Bloch regions, {\it Comput. Methods Funct. Theory} (1994) (Penang), 89--100, Ser. Approx. Decompos., {\bf 5}, {\it World Sci. Publ., River Edge, NJ}, 1995.
	
	\bibitem{Bonk-J-Analyze-Math-1996} {\sc M. Bonk, D. Minda,} and {\sc H. Yanagihara}, Distortion theorems for locally univalent Bloch functions, {\it J. Anal. Math.} {\bf 69} (1996), 73--95.
	
	
	
	\bibitem{chen-PAMS-2000} {\sc H. Chen, P. Gauthier} and {\sc W. Hengartner}, Bloch constants for planar harmonic mappings, {\it Proc. Amer. Math. Soc.} {\bf 128} (2000), 3231--3240.
	
	\bibitem{chen-appl-math-comp-2009} {\sc S. Chen, S. Ponnusamy} and {\sc X. Wang},  Landau’s theorem for certain biharmonic mappings, {\it Appl. Math. Comput.} {\bf 208} (2009), 427-433.
	
	\bibitem{chen-2011} {\sc S. Chen, S. Ponnusamy} and {\sc X. Wang}, Landau's theorem and Marden constant for harmonic $\nu$-Bloch mappings, {\it Bull. Aust. Math. Soc.} {\bf 84} (2011), 19--32.
	
	\bibitem{chen-2014-JMMA} {\sc J. Chen, A. Rasila} and {\sc X. Wang}, Landau’s theorem for polyharmonic mappings, {\it  J. Math. Anal. Appl.} {\bf 409} (2014), 934--945.
	
	\bibitem{colona-1987} {\sc F. Colonna}, The Bloch constant of bounded analytic functions, {\it J. London Math. Soc.} {\bf 36} (1987), 95--101.
	
	\bibitem{colona-1989} {\sc F. Colonna}, The Bloch constant of bounded harmonic mappings, {\it Indiana Univ. Math. J.} {\bf 38} (1989), 829--840.
	
	\bibitem{defant-2003} {\sc A. Defant}, {\sc D. Garc\'{i}a}, and {\sc M. Maestre},  Bohr's power series theorem and local Banach space theory, {\it J. Reine Angew. Math.} \textbf{557} (2003), 173–197.
	
	\bibitem{defant-2006} {\sc A. Defant} and {\sc L. Frerick}, A logarithmic lower bound for multi-dimenional bohr radii, {\it Israel J. Math.} {\bf 152} (2006), 17--28.
	
	\bibitem{defant-JRAM-2011} {\sc A. Defant} and {\sc L. Frerick}, The Bohr radius of the unit ball of $l^{n}_{p}$, {\it J. Reine Angew. Math.} {\bf 660} (2011), 131--147.
	
	\bibitem{defant-2011} {\sc A. Defant, L. Frerick, J. Ortega-Cerd${\rm \grave{A}}$, M. Ouna${\rm \ddot{I}}$es}, and {\sc K. Seip}, The Bohnenblust-Hille inequality for homogeneous polynomils in hypercontractive, {\it Ann. of Math.} {\bf 174} (2011), 512--517.
	
	
	\bibitem{Dixon & BLMS & 1995} {\sc P. G. Dixon}, Banach algebras satisfying the non-unital von Neumann inequality, {\it Bull. London Math. Soc.} \textbf{27} (4) (1995), 359--362.
	
	\bibitem{efraimidis-2017} {\sc I. Efraimidis, J. Gaona, R. Hern${\rm \acute{A}}$ndez} and {\sc O. Venegas}, On harmonic Bloch-type mappings, {\it Complex Var. Elliptic Equ.} {\bf 62} (2017), 1081--1092.
	 
	
	
	
	\bibitem{Evd-Ponn-Rasi-2020} {\sc S. Evdoridis}, {\sc S. Ponnusamy} and {A. Rasila}, Improved Bohr's inequality for shifted disks, {\it Results Math.} {\bf 76:14} (2021), 15 pages.
	
	\bibitem{Fernandez-1984} {\sc J. L. Fern\'{a}ndez}, On the coefficients of Bloch functions, {\it J. London Math. Soc.} {\bf 29} (1984), 94--102.
	
	
	
	\bibitem{Four-Rusc-2010} {\sc R. Fournier} and {\sc St. Ruscheweyh}, On the Bohr radius for simply connected domains, \textit{Centre de Recherches Math$ \acute{e} $matiques CRM Proceedings and Lecture Notes}, Vol. \textbf{51} (2010), 165--171.
	
	
	
	
	
	
	
	\bibitem{Gnuschke-Hauschild-1986} {\sc D. Gnuschke-Hauschild} and {\sc Ch. Pommerenke}, On Bloch functions and gap series, {\it J. Reine Angew. Math.} \textbf{367} (1986), 172–186.
	
    
	
	\bibitem{Ismagilov-2020} {\sc A. Ismagilov, I. R. Kayumov,} and {\sc S. Ponnusamy}, Sharp Bohr type inequality, 
	{\it J. Math. Anal. Appl.}  {\bf 489} (2020), 124147.
	
	\bibitem{kalaj-2014} {\sc D. Kalaj}, {\sc S. Ponnusamy}, and {\sc M. Vuorinen}, Radius of close-to-convexity and fully starlikeness of harmonic mappings, {\it  Complex Var. Elliptic Equ.} {\bf 44} (2014), 685--692.
	
	
	
	
	
	
	\bibitem{Kayumov-Ponnusamy-2018-b} {\sc I. R. Kayumov} and {\sc S. Ponnusamy}, Bohr's inequalities for the analytic functions with lacunary series and harmonic functions, 
	{\it J. Math. Anal. Appl.}  {\bf 465} (2018), 857--871.
	
	\bibitem{Kay & Pon & AASFM & 2019} {\sc I. R. Kayumov} and {\sc S. Ponnusamy}, On a powered Bohr inequality, \textit{Ann. Acad. Sci. Fenn. Ser. A}, \textbf{44}(2019), 301--310.
	
	
	
	
	\bibitem{Kayumov-Ponnuswamy-MN-2018} {\sc I. R. Kayumov}, {\sc S. Ponnusamy} and {\sc N. Shakirov}, Bohr radius for locally univalent harmonic mappings, {\it Math. Nachr.} {\bf 291} (2018), 1757--1768.
	
	
	
	
	
	
	
	
	
	\bibitem{Liu-Results-Math-2018} {\sc G. Liu} and {\sc S. Ponnusamy}, On Harmonic $\nu$-Bloch and $\nu$-Bloch-type mappings, {\it Results Math.} {\bf 73:90} (2018), 21 pages.
	
	
	
	\bibitem{Liu-Ponnusamy-BMMS-2019} {\sc Z. H. Liu} and {\sc S. Ponnusamy}, Bohr radius for subordination and $ K $-quasiconformal harmonic mappings, {\it Bull. Malys. Math. Sci. Soc.} \textbf{42} (2019), 2151--2168.
	

\bibitem{M-S-Liu-Landau-2009} {\sc M. S. Liu}, Landau’s theorems for planar harmonic mappings, {\it  Comput. Math. Appl.} {\bf 57} (2009), 1142--1146.


\bibitem{M-S-Liu-bloch-2009} {\sc M. S. Liu}, Estimates on Bloch constants for planar harmonic mappings, {\it Sci. China Ser. A Math.} {\bf 52} (2009), 87--93.
	
	\bibitem{Liu-Pon-PAMS-2020} {\sc M. S. Liu} and {\sc S. Ponnusamy}, Multidimensional analogues of refined Bohr's inequality, {\it Proc. Amer. Math. Soc.} {\bf 149} (2021), 2133--2146.
	
	
	
	
	
	
		\bibitem{paulsen-2002} {\sc V. I. Paulsen, G. Popescu} and {\sc D. Singh}, On Bohr's inequality, {\it Proc. Lond. Math. Soc.} s3-85 (2002), 493--512.
	
	
	
	
	\bibitem{popescu-2019} {\sc G. Popescu}, Bohr inequalities for free holomorphic functions on polyballs, {\it Adv. Math.} {\bf 347} (2019), 1002-1053.
	
	
	
	
	
	
	\bibitem{Rohde-1993} {\sc S. Rohde}, The boundary behavior of Bloch functions, {\it J. London Math. Soc.} {\bf 48} (1993), 488–-499. 
	
	
	
	\bibitem{Ruscheweyh-1985} {\sc St. Ruscheweyh}, Two remarks on bounded analytic functions, \emph{Serdica} \textbf{11}(1) (1985), 731--732.
	
	
	
	
	
	
	
	
	
	
	
	\bibitem{Yanigaha-BLMS-1994} {\sc H. Yanagihara}, Sharp distortion estimate for locally schlicht Bloch functions, {\it Bull. London Math. Soc.} {\bf 26} (1994), 539--542.
	
	\bibitem{zhu-CAOT-2015} {\sc J. F. Zhu}: Landau theorem for planar harmonic mappings, {\it Complex Anal. Oper. Theory} {\bf 9} (2015), 1819--1826.
	
	\bibitem{zhu-2016} {\sc J. Zhu}, Coefficients estimate for harmonic $\nu$-Bloch-mappings and harmonic $K$-quasiconformal mappings, {\it Bull. Malays. Math. Sci. Soc.} {\bf 39} (2016), 349--358.
\end{thebibliography}
\end{document}